\documentclass[a4paper,10pt]{article} 
\usepackage[utf8]{inputenc} 
\usepackage[T1]{fontenc}   
\usepackage[english]{babel} 

\usepackage[backend=biber,style=ieee-alphabetic,isbn=false,doi=false,url=false, clearlang=true,sorting=nyt]{biblatex}
\DeclareFieldFormat[article]{title}{\textit{#1}}
\DeclareFieldFormat[article]{date}{(#1)}
\DeclareFieldFormat[article]{journaltitle}{\*{#1}}
\DeclareFieldFormat[article]{volume}{\textbf{#1}}
\DeclareFieldFormat[book]{title}{\*{#1}}
\DeclareFieldFormat[book]{number}{vol. {#1}}
\DeclareFieldFormat[incollection]{title}{\textit{#1}}
\DeclareFieldFormat[incollection]{booktitle}{\*{#1}}
\AtEveryBibitem{\clearlist{language}}

  \renewbibmacro*{journal+issuetitle}{%
  \usebibmacro{journal}%
  \setunit*{\addspace}%
  \iffieldundef{series}
    {}
    {\newunit
     \printfield{series}%
     \setunit{\adddot\space}}
      \printfield{volume}%
      \setunit*{\addspace}
      \usebibmacro{date}
        \setunit{\addcomma\space}
      \printfield{number}
        \setunit{\addcomma\space}
      \usebibmacro{issue}
  \newunit}

\DeclareBibliographyDriver{article}{%
  \usebibmacro{bibindex}%
  \usebibmacro{begentry}%
  \usebibmacro{author/translator+others}%
  \setunit{\printdelim{nametitledelim}}\newblock
  \usebibmacro{title}%
  \newunit
  \printlist{language}%
  \newunit\newblock
  \usebibmacro{byauthor}%
  \newunit\newblock
  \usebibmacro{bytranslator+others}%
  \newunit\newblock
  \printfield{version}%
  \newunit\newblock
  \usebibmacro{journal+issuetitle}%
  \newunit
  \usebibmacro{byeditor+others}%
  \newunit
  \usebibmacro{note+pages}%
  \newunit\newblock
  \iftoggle{bbx:isbn}
    {\printfield{issn}}
    {}%
  \newunit\newblock
  \usebibmacro{doi+eprint+url}%
  \newunit\newblock
  \usebibmacro{addendum+pubstate}%
  \setunit{\bibpagerefpunct}\newblock
  \usebibmacro{pageref}%
  \newunit\newblock
  \iftoggle{bbx:related}
    {\usebibmacro{related:init}%
     \usebibmacro{related}}
    {}%
  \usebibmacro{finentry}}

  \DeclareBibliographyDriver{book}{%
  \usebibmacro{bibindex}%
  \usebibmacro{begentry}%
  \usebibmacro{author/editor+others/translator+others}%
  \setunit{\printdelim{nametitledelim}}\newblock
  \usebibmacro{maintitle+title}%
  \newunit
  \printlist{language}%
  \newunit\newblock
  \usebibmacro{byauthor}%
  \newunit\newblock
  \usebibmacro{byeditor+others}%
  \newunit\newblock
  \printfield{edition}%
  \newunit
  \iffieldundef{maintitle}
    {\printfield{volume}%
     \printfield{part}}
    {}%
  \newunit
  \printfield{volumes}%
  \newunit\newblock
  \usebibmacro{series+number}%
  \newunit\newblock
  \printfield{note}%
  \newunit\newblock
  \usebibmacro{publisher+location+date}%
  \newunit\newblock
  \usebibmacro{chapter+pages}%
  \newunit
  \printfield{pagetotal}%
  \newunit\newblock
  \iftoggle{bbx:isbn}
    {\printfield{isbn}}
    {}%
  \newunit\newblock
  \usebibmacro{doi+eprint+url}%
  \newunit\newblock
  \usebibmacro{addendum+pubstate}%
  \setunit{\bibpagerefpunct}\newblock
  \usebibmacro{pageref}%
  \newunit\newblock
  \iftoggle{bbx:related}
    {\usebibmacro{related:init}%
     \usebibmacro{related}}
    {}%
  \usebibmacro{finentry}}

\renewbibmacro*{publisher+location+date}{%
  \printlist{location}%
  \iflistundef{publisher}
    {\setunit*{\addcomma\space}}
    {\setunit*{\addcolon\space}}%
  \printlist{publisher}%
  \setunit*{\addcomma\space}%
  \usebibmacro{date}
  \adddot%
  \newunit}

\DeclareBibliographyDriver{incollection}{%
  \usebibmacro{bibindex}%
  \usebibmacro{begentry}%
  \usebibmacro{author/translator+others}%
  \setunit{\printdelim{nametitledelim}}\newblock
  \usebibmacro{title}%
  \newunit
  \printlist{language}%
  \newunit\newblock
  \usebibmacro{byauthor}%
  \newunit\newblock
  \usebibmacro{in:}%
  \usebibmacro{maintitle+booktitle}%
  \newunit\newblock
  \usebibmacro{byeditor+others}%
  \newunit\newblock
  \printfield{edition}%
  \newunit
  \iffieldundef{maintitle}
    {\printfield{volume}%
     \printfield{part}}
    {}%
  \newunit
  \printfield{volumes}%
  \newunit\newblock
  \usebibmacro{series+number}%
  \newunit\newblock
  \printfield{note}%
  \newunit\newblock
    \usebibmacro{chapter+pages}%
  \newunit\newblock
  \iftoggle{bbx:isbn}
    {\printfield{isbn}}
    {}%
  \newunit\newblock
  \usebibmacro{publisher+location+date}%
  \newunit\newblock
  \usebibmacro{doi+eprint+url}%
  \newunit\newblock
  \usebibmacro{addendum+pubstate}%
  \setunit{\bibpagerefpunct}\newblock
  \usebibmacro{pageref}%
  \newunit\newblock
  \iftoggle{bbx:related}
    {\usebibmacro{related:init}%
     \usebibmacro{related}}
    {}%
  \usebibmacro{finentry}}

\usepackage{mathtools,amssymb,amsthm}
\usepackage{mathrsfs,mathdots}
\usepackage{stmaryrd}
\usepackage{dsfont}
\usepackage{a4wide}
\usepackage{MnSymbol}
\usepackage{pst-all,pstricks,pstricks-add,pst-plot,pst-eucl}
\usepackage[all]{xy}
\usepackage{faktor}
\usepackage{xcolor}
\usepackage{ytableau}
\usepackage{braket}
\usepackage{hyperref} 
\allowdisplaybreaks
\newtheorem{theo}{Theorem}[section]
\newtheorem{lem}[theo]{Lemma}
\newtheorem{prop}[theo]{Proposition}
\newtheorem{cor}[theo]{Corollary}
\newtheorem{defi}[theo]{Definition}
\newtheorem{rqe}[theo]{Remark}
\newtheorem{exemple}[theo]{Example}

\newcommand{\conj}[2]{{}^{#1}#2}

\def\a{\alpha}
\def\e{\varepsilon}
\def\d{\delta}
\def\s{\sigma}
\def\w{\omega}
\def\g{\mathfrak{g}}

\def\B{\mathcal{B}}
\def\L{\mathcal{L}}

\def\s{\sigma}
\def\k{\mathds{k}}

\def\F{\mathcal{F}}
\def\G{\mathcal{G}}
\def\H{\mathcal{H}}
\def\T{\mathcal{T}}
\def\Z{\mathcal{Z}}
\def\P{\mathcal{P}}
\def\L{\mathcal{L}}
\def\O{\mathcal{O}}

\def\W{\mathcal{W}}
\def\X{\mathcal{X}}

\def\ZZ{\mathbb{Z}}

\newcommand\smallw{
  \mathchoice
    {{\scriptstyle\mathcal{W}}}
    {{\scriptstyle\mathcal{W}}}
    {{\scriptscriptstyle\mathcal{W}}}
    {\scalebox{.5}{$\scriptscriptstyle\mathcal{W}$}}
  }
  
  \DeclareFontFamily{U}{mathx}{\hyphenchar\font45}
\DeclareFontShape{U}{mathx}{m}{n}{
      <5> <6> <7> <8> <9> <10>
      <10.95> <12> <14.4> <17.28> <20.74> <24.88>
      mathx10
      }{}
\DeclareSymbolFont{mathx}{U}{mathx}{m}{n}
\DeclareFontSubstitution{U}{mathx}{m}{n}
\DeclareMathAccent{\widecheck}{0}{mathx}{"71}
\DeclareMathAccent{\wideparen}{0}{mathx}{"75}

\begin{document}

\pdfbookmark[0]{Page de garde}{garde}
\thispagestyle{empty}

\title{Orbit closures in flag varieties for the centralizer of an order-two nilpotent element: normality and resolutions for types A, B, D. }
\author{Simon Jacques}
\maketitle
\abstract{Let $G$ be a reductive algebraic group in classical types A, B, D. Let $e$ be an element of the Lie algebra of $G$ with $Z\subset G$ its centralizer for the adjoint action. We assume that $e$ identifies with a nilpotent matrix of order two, which guarantees that the number of $Z$-orbits in the flag variety of $G$ is finite. For types B and D in characteristic two, we also assume that the image of $e$ is totally isotropic. We show that the closure $Y$ of such an orbit is normal. We also prove that $Y$ is Cohen-Macaulay with rational singularities provided that the base field is of characteristic zero, and that Cohen-Macaulayness holds in any characteristic for type A. We exhibit a rational and birational morphism onto $Y$ involving Schubert varieties. Our work generalizes a result by N. Perrin and E. Smirnov on the Springer fibers.}

\section*{Introduction}

\paragraph*{1.}
Let $\k$ be an algebraically closed field and let $G$ be a reductive connected algebraic group over $\k$ with $B$ a Borel subgroup. Let $e$ be a nilpotent element of the Lie algebra $\g$ of $G$, and let $Z$ be its centralizer in $G$ for the adjoint action. When the number of $Z$-orbits in the flag variety $G/B$ is finite, their closures are of particular interest. They include in this case the irreducible components of the so called Springer fiber over $e$. It is the fiber of $e$ under the proper birational morphism $$\overset{\sim}{\mathcal{N}}\to \mathcal{N},$$ called the Springer resolution, which is the projection onto the nilpotent cone $\mathcal{N}\subset \g$ from the smooth variety $\overset{\sim}{\mathcal{N}}:=\set{(x,gB)\in  \mathcal{N}\times G/B| Ad\ g^{-1}\cdot x\in \mathfrak{b}}$, $\mathfrak{b}$ denoting the Lie algebra of $B$. The Springer fibers are of main interest in representation theory (see the seminal work of T. A. Springer \cite{Springer1976}, their link with the orbital varieties \cite{Spaltenstein1982} and the Steinberg variety \cite{Steinberg1976}\footnote{Actually, the latter references deal with unipotent elements instead of nilpotent ones, regarding the Springer fibers as the variety of Borel subgroups containing a given unipotent element. However, recall that when $G$ is the general linear group or is almost simple and simply connected and the characteristic of $\k$ is good, the unipotent variety in the group $G$ and the nilpotent cone in its Lie algebra $\g$ can be identified with a $G$-equivariant isomorphism (see for example \cite[Theorem 3.1]{Springer1976} for an original but weaker statement and \cite[Theorem 6.20]{Humphreys1995} and \cite[Corollary 9.3.3]{BardsleyRichardson1985} for this more general one), so that the two notions of Springer fibers exactly match. }). They are connected and equidimensional (see for example \cite{Spaltenstein1982}) and their irreducible components have been the subject of numerous studies. For the classical cases and $char(\k) \neq 2$, N. Spaltenstein (type A, \cite{Spaltenstein1982}) and M. van Leeuwen (types B, C, D, \cite{Leeuwen1989}) showed that they are parameterized by standard and domino tableaux, whose shapes are given by Young diagrams relative to the nilpotent orbit in question. Subsequent studies of their singularities have often been based on these shapes and have mainly produced results for $G$ the general linear group and $\k$ the field of complex numbers. For an example, F. Fung showed in \cite{Fung} that they are all smooth, in the so-called hook and two-line cases. A. Melnikov and L. Fresse gave a necessary and sufficient condition for this global smoothness in \cite{FresseMelnikov2010} while they gave a criterion for individual smoothness in \cite{Fresse2009,FresseMelnikov2010}, under the additional assumption of being in the \textit{two column case}. This is the first case where singularities appear. It also implies that the order of nilpotency of $ad\ e$ is less than or equal to $3$, which is a condition ensuring, after the work of Panyushev \cite{Panyushev1994}, that the number of $Z$-orbits is finite (in fact, this implication is established for $char(\k)=0$, but for our types A, B, D, it is still valid for the other characteristics, see Propositions \ref{prop:type A T fixed point} and \ref{prop:types BD T fixed point}).

\paragraph{2.} The \textit{two-column case} is assumed in the article \cite{PerrinSmirnov} by N. Perrin and E. Smirnov. For type A and $char(\k) \neq 2$, they present rational resolutions of the components and show that they are normal and Cohen-Macaulay. They also give arguments for the same results in type D, but there is a gap in their proof of normality and the Cohen-Macaulay property, due to the non-algebraicity of a certain map (see Appendix B for details and a counter-example). Nevertheless, their proof of the existence of a rational birational morphism onto the component is still valid for this type. Our work is mainly inspired by the latter, generalizing it in several directions. Retaining the assumption of the two-column case, we also prove normality and rationality, but for the much broader class of $Z$-orbit closures. For example, if $G=Gl_{n\k}$ is the general linear group and $r$ represents the rank of $e$ considered as a nilpotent matrix of order two, then we can deduce from our Proposition \ref{prop:type A T fixed point} and the hook-length formula that the number of $Z$-orbits is $(n-r+1)(n-r)\dots (n-2r+2)$ times the number of irreducible components. In addition, we consider the three types A, B, D and we also deal with the case $char(\k)=2$ (with precautions regarding the nilpotent orbit considered, see below).

\paragraph{3.} Let us now state our main results. We assume that $\k$ is of arbitrary characteristic and we fix an integer $n$. Let ${O_{n}}_{\k}$ be the group over $\k$ whose closed points are the invertible $n\times n$ matrices preserving the quadratic form
\begin{equation}\label{eq:forme quadratique générale}
\sum_{k=1}^{ \lfloor (n+1)/2 \rfloor }Y_{k}Y_{n-k+1}.
\end{equation}
For even $n$, let us denote by $\vartriangle_{n}$ the Dickson invariant, as defined for example in \cite[IV, §5]{Knus}. This is the regular function on ${O_{n}}_{\k}$ satisfying $det_{n}=1+2\vartriangle_{n}$ where $det_{n}$ is the restriction of the determinant to ${O_{n}}_{\k}$ (see Section \ref{sec:matrix model BCD} for details). We then define the special orthogonal group ${SO_{n}}_{\k}$ as the zero locus of $\vartriangle_{n}$ if $n$ is even and as that of $det_n-1$ if $n$ is odd. Without these precautions, note that it fails to be connected and semi-simple of type B$_n$ (odd $n$) or D$_n$ (even $n$) in the case $char(\k)=2$ (see for example \cite{Hesselink1979}, \cite[Appendix C]{ConradReductiveGroupScheme} and Section \ref{sec:comments matrix models}).

Assume now that $G$ is the general linear group ${Gl_{n}}_{\k}$ or the special orthogonal group ${SO_{n}}_{\k}$. We also assume that the nilpotent element $e$ is identified with a nilpotent matrix of order two, which means that we are in the two-column case.  If $char(\k)=2$ and $G={SO_{n}}_{\k}$, we make the additional assumption that the image of $e$ is totally isotropic.

Recall that a proper morphism $f:X\to Y$ of locally noetherian schemes is called \textit{rational} if $\O_Y\simeq f_*\O_X$ and $R^i f_* \O_X=0$ for $i>0$. When the schemes are irreducible varieties with $X$ smooth, such a rational morphism $f$ is said to be a \textit{rational resolution} if it is also birational with $R^i f_* \w_{X}=0$ for $i>0$, where $\w_{X}$ denotes the canonical bundle of $X$. If $char(\k)=0$, two rational resolutions can be dominated by a third, so being the target of a rational resolution leads to the intrinsic notion of \textit{having rational singularities}. We prove:

\begin{theo}\label{theo:0 normality, Cohen-Macaulay etc.}
The $Z$-orbit closures in the flag variety of $G$ are normal. In characteristic zero, they are Cohen-Macaulay with rational singularities. In any other characteristic, they remain Cohen-Macaulay in type A. 
\end{theo}

This theorem is based on two results. The first is the construction of an explicit birational morphism using matrix models and involving Schubert varieties. It ensures the existence of a Borel subgroup $B$ of $G$, containing a maximal torus $T$, and of a closed reductive subgroup $H$ of $G$ equipped with a retraction $\varpi:Z\to H$, having $B_H:=B\cap H$ as Borel subgroup and $T_H:=T\cap H$ as maximal torus, so that

\begin{theo}\label{theo:1 birational rational}
For any $Z$-orbit closure $Y$ in $G/B$, there exists $w$ in the Weyl group of $G$ such that $Y=\overline{HB\cdot wB}=\overline{Z\cdot wB}$ and 
\begin{equation}\label{eq:morphisme birationnel}
H\times^{B_H}\overline{B\cdot wB}\to Y,\ [h,gB]\mapsto hgB
\end{equation} is rational, birational, $Z$-equivariant, with a $Z$-action on $H\times^{B_H}\overline{B\cdot wB}$ defined by $z\cdot [h,gB]=[\varpi(z)h,h^{-1}\varpi(z)^{-1}zhgB]$.
\end{theo}

The second result is valid in a more general context, where we assume only that $G$ is a connected reductive group over $\k$ and $H$ a closed connected reductive subgroup of $G$. We make the same assumptions as before about $T$, $B$, $T_H$, $B_H$ and fix any $w$ in the Weyl group of $G$.  We denote by $\rho_G$ the half-sum of positive roots and, for any dominant character $\lambda$, by $V_G\left(\lambda\right)$ the dual Weyl $G$-module with lowest weight $-\lambda$. Let $\rho_H$ and $V_H\left(\lambda\right)$ also be the corresponding objects for $H$. We refer to Section \ref{sec:normality rationality} for details of the notation and a stronger result that also deals with the vanishing of the canonical bundle.

\begin{theo} \label{theo:2 rationality normality general argument}
\begin{enumerate} Let us assume that 
\item[(i)] the morphism $\pi\colon H\times ^{B_H} \overline{B\cdot wB}\to \overline{HB\cdot wB}, [h,gB]\mapsto hgB$ is birational,
\item[(ii)] the character $2\rho_H-\rho_{G\vert T_H}$ is dominant,
\item[(iii)] $char(\k)=0$ or 
\item[(iii)'] $char(\k)=p>0$ and the restriction $V_G\left((p-1)\rho_G\right)\to V_H\left((p-1)\rho_{G\vert T_H}\right)$ is surjective.
\end{enumerate}
Then $\overline{HB\cdot wB}$ is normal and $\pi$ is rational. 
\end{theo}

\begin{rqe}
It remains an open question whether the Cohen-Macaulay property is valid for types B and D and whether this property, rationality and normality are valid for type C and the exceptional types.
\end{rqe}

\begin{rqe}
If we take $H=T$ in Theorem \ref{theo:2 rationality normality general argument}, we find the well-known result on the normality of Schubert varieties. In fact, in this case, the sequence of arguments used in the proof coincides with that of M. Brion and S. Kumar in \cite{BrionKumar}. 
\end{rqe}

\begin{rqe}
In types B and D, if $char(\k)\neq 2$, the assumption of being in the two-column case implies that the image of $e$ is totally isotropic. However, this is not the case if $char(\k)=2$, as can be seen in type D by taking $e:=\begin{psmallmatrix}0 & 1 & 1 & 0\\
1 & 0 & 0 & 1\\
1 & 0 & 0 & 1\\
0 & 1 & 1 & 0
\end{psmallmatrix}$ which is a matrix of nilpotency order two as $e':=\begin{psmallmatrix}0 & 0 & 1 & 0\\
0 & 0 & 0 & 1\\
0 & 0 & 0 & 0\\
0 & 0 & 0 & 0
\end{psmallmatrix}$. Note that the dimensions of the centralizers of $e$ and $e'$ differ (they are $2$ and $4$ respectively, see \cite[Theorem 4.5]{Hesselink1979}), so we can see that the rank of nilpotent elements does not suffice to characterize nilpotent orbits in this case. Our assumption about the image of $e$ is therefore necessary to work in our matrix model, and then apply our reasoning (which depends crucially on the dimension of $Z$). It also turns out to be sufficient (see Section \ref{sec:application through matrix models} and the result of \cite[Theorem 3.8]{Hesselink1979}).
\end{rqe}

\begin{rqe}\label{rqe:type C}
Let us assume that $char(\k)\neq 2$, that the rank of $e$ is odd and that the $Z$-orbit in question has a $T$-fixed point (which is not superfluous, since there are orbits without such points, see Proposition \ref{prop:types BD T fixed point}). Finally, let us replace $Z$ by its neutral component $Z^0$. Theorem \ref{theo:1 birational rational} (with the exception of rationality) is then again valid for $G={Sp_n}_{\k}$ the symplectic group over $\k$ (see the matrix models in Section \ref{sec:matrix models} for a concrete description of this group). However, our proof of Theorem \ref{theo:0 normality, Cohen-Macaulay etc.} cannot work because of the non-dominance of the character involved in Theorem \ref{theo:2 rationality normality general argument} (see Remark \ref{rqe:type C caractère non dominant}).
\end{rqe}

\paragraph{4.} To prove these results, we take up many of Perrin and Smirnov's arguments in \cite{PerrinSmirnov} while developing new techniques. Our birational morphism (\ref{eq:morphisme birationnel}) onto a $Z$-orbit closure is quite analogous to Perrin and Smirnov's morphism targeting an irreducible component of the Springer fiber (see (\ref{eq:PerrinSmirnov morphism}) in Appendix B); in particular, each morphism involves a Schubert variety. However, our approach is different. Indeed, Perrin and Smirnov present an explicit description and a direct proof of birationality, depending on the type A or D and adapted to the particular framework of the irreducible components of the Springer fiber. Moreover, the presence of a Schubert variety is a simple consequence of their construction; it is detected in their source variety as the fiber of a certain projection. For us, in contrast, the Schubert variety is a starting point for producing more general resolutions, enabling us to deal simultaneously
with types A, B, C, D (see Lemma \ref{lem:équivalence birationnelle cas général} with the additional assumptions concerning type C, as indicated in the Remark \ref{rqe:type C}.), and recovering for type A, modulo some precautions, the previous construction (see (\ref{eq:identification birational Perrin Smirnov}) in Appendix B). To do this, we use results on parabolic induction (see \cite{ChaputFresse} by P-E. Chaput, L. Fresse, T. Gobet) and symmetric subgroups (see the work of R. W. Richardson and T. A. Springer \cite{RichardsonSpringer1990}). Nevertheless, this general construction remains highly technical, since it is designed to be applied to a very specific matrix model.

Based on the previous birationality statement, our proof of normality, rationality and Cohen-Macaulay property closely follows that of Perrin and Smirnov, except for one point. We use the same reference \cite{HeThomsen} of X. He and J-F. Thomsen to establish a Frobenius splitting\footnote{Note that our approach is a little simpler since we use the main theorem of \cite{HeThomsen} whereas Perrin and Smirnov, in order to obtain a more precise splitting, combine several results from this reference.} in order to obtain the surjectivity of a certain restriction of sections (see (\ref{eq:surjectivity of restriction}) in Section \ref{sec:normality and rationality}), and finally the desired results using the same inductive argument (Proposition \ref{prop:PerrinSmirnov induction}). Moreover, our calculation for proving the Cohen-Macaulay property also follows their (with slight adaptations for types B, C, D, see Proposition \ref{prop:canonical sheaf formula}). 

The main difference lies in the arguments for extending the consequences of the Frobenius splitting from the positive characteristic to all characteristics. One way of doing this is to \textit{realize} the desired varieties in mixed characteristics. For Perrin, Smirnov and us, this will mean producing flat schemes of finite presentation on sufficiently large bases, which are such that the collections of their geometric fibers account for different incarnations of the starting varieties, in positive and zero characteristics. But by virtue of the equations that define them, the Springer fibers and their irreducible components can be easily realized over $\mathbb{Z}$, and Perrin and Smirnov did not need to resort to more complicated arguments.

As for us, we must realize all varieties of the form $\overline{HB\cdot wB}$. This was done for $H=T$ (Schubert varieties) by V-B. Mehta, A. Ramanathan in \cite{MehtaRamanathan1985}. But the general case requires systematically dealing with a problem of scheme-theoretic image formation under non-flat base changes. By abandoning $\mathrm{\mathrm{Spec}}\ \mathbb{Z}$ for smaller bases $\mathrm{\mathrm{Spec}}\ A$ where $A$ is an integral algebra of finite type over $\mathbb{Z}$, we present a solution for scheme-theoretic image realizations under some assumptions (such that properness and having integral geometric fibers, see Theorem \ref{theo:relèvement image}). This leads to fairly general realizations of closures, in a framework of actions of $\k$-groups on proper $\k$-schemes of finite type (Corollary \ref{cor:réalisation H1H2...HnZ barre}) including the case of our $\overline{HB\cdot wB}$ varieties (see Lemma \ref{lem:realisation image intermediaire} and its proof in Appendix A).

\paragraph{5.} If the number of $Z$-orbits is finite, another interesting fact is that the $G/Z$-variety is spherical (i.e., it has a finite number of $B$-orbits). Transposed in terms of $B$-orbit closures in $G/Z$, our study then fits into the theory of spherical varieties, as originally developed by D. Luna in \cite{Luna2001} and then by F. Knop in \cite{Knop2014}. In \cite{Brion2003}, M. Brion presents for such $B$-orbit closures a powerful criterion for normality and Cohen-Macaulay property, called the multiplicity free criterion. Applying it to our situation remains an open and interesting problem.

\paragraph{6.} Our article is organized as follows. We first introduce matrix models for type A on the one hand and types B, C, D on the other (Section \ref{sec:matrix models}). We then prove Theorem \ref{theo:1 birational rational} with the exception of the rationality hypothesis, by applying two technical lemmas on matrix models (Section \ref{sec:birationality}). With Theorem \ref{theo:theorem 2 plus} in Section \ref{sec:normality rationality} we prove (a strong version of) Theorem \ref{theo:2 rationality normality general argument}, using our realization result presented in Appendix A. Some additional checks on matrix models allow us to combine theorems \ref{theo:1 birational rational} and \ref{theo:theorem 2 plus} and conclude in Section \ref{sec:conclusion} to our main result (Theorem \ref{theo:0 normality, Cohen-Macaulay etc.}). In Appendix B, we present some comments on the comparison with Perrin and Smirnov's paper.

\paragraph{7.} Throughout this article, the varieties are reduced schemes of finite type over an algebraically closed field, and the (linear) algebraic groups are affine group schemes which are varieties. In the usual way, we often only consider closed points when working on varieties, and in this context, taking the $\cap$ intersection means we are considering the reduced structure on the scheme-theoretic intersection (fiber product). We refer to \cite[Exposés XIX to XXVI]{SGA3XIXaXXVI} for the definitions of algebraic group notions (Borel subgroup, maximal Torus, being reductive, ...) in the context of group schemes. 

\paragraph{Acknowledgment}
This article comes from a PhD thesis made under the direction of P-E. Chaput and L. Fresse. It could not have been done without them. They inspired most of its ideas and took a lot of time to proofread it; it was moreover a true pleasure to work under their guidance and I am in a considerable debt to them. I would express my gratitude to N. Perrin and E. Smirnov whose work is the constant reference of this article, and with whom we had several fruitful discussions. I would also thank A. Genestier for its support and M. Romagny for his advice and help. S. Cupit-Foutou, A. Moreau and K. Česnavičius were the members of my thesis jury with E. Smirnov, A. Genestier and M. Romagny, and I thank all of them for their questions and reports on my work. Finally, I would like to thank E. Zabeth, A. Lacabanne and M. Metodiev for their generous help. 

\begin{center}\textit{Apart from Appendix A, $\mathds{k}$ will denote an algebraically closed field.}\end{center}
{\tableofcontents}

{\vspace{\stretch{3}}}

\section{Matrix models}\label{sec:matrix models}

\subsection{General notation}

We first introduce the following notation and conventions.
\begin{enumerate}
\item[-] For any group $G$ (group scheme or algebraic group) we will say that $(T,B)$ is a Killing pair if $B$ is a Borel subgroup of $G$ containing a maximal torus $T$. If we fix such a pair and the group is assumed to be reductive, we denote by $\rho_G$ the sum of all fundamental weights, that is the half sum of positive roots.  If $P\supset B$ is a parabolic subgroup of $G$, we denote by $W_P$ its Weyl group relative to $T$ and by $W^P$ the system of minimum-length representatives of the quotient $W/W_P$.
\item[-] If there is no ambiguity, we set $\overline{i}:=m-i+1$.  We define $\check{\sigma}$ as the permutation $i\mapsto \overline{\sigma (\overline{i})}$ of $\set{1,...,m}$. We denote by $\ell_m (\s)$ the number of inversions of $\s$, that is its length in $\mathfrak{S}_m$: $\ell_m(\s):=\#\set{1\leq i<j\leq m| \s(i)>\s(j)}$. For $\varsigma \in \mathfrak{S}_q$, $m\leq q$, $0\leq k\leq q-m$ we say that $\sigma\in \mathfrak{S}_m$ is \textit{the induced permutation of} $\varsigma \in \mathfrak{S}_q$ \textit{on} $\set{k+1,...,k+m}\subset \set{1,...,q}$ if the restriction of $\varsigma\begin{psmallmatrix}
I_{k} & 0 & 0\\
0 & \s^{-1} & 0\\
0 & 0 & I_{q-(k+m)}
\end{psmallmatrix}$ on $\set{k+1,...,k+m}$ is increasing. 
\item[-] For a matrix $M$, we denote by ${}^{\d}M$ the symmetric transform of $M$ along the antidiagonal (namely, exchanging coefficients $(i,j)$ and $(\overline{j},\overline{i})$). For any ring $R$, we often identify the permutation group $\mathfrak{S}_n$ with its image in $Gl_n(R)$, thanks to the action on $R^n$ given by $\sigma \cdot (v_i)_i=(v_{\sigma(i)})_i$. In this point of view, we have $\check{\s}={}^{\delta}\sigma ^{-1}$ for any permutation $\s$. If $m\in \mathbb{N}$ and $\epsilon=\pm 1$, we denote $I_{\epsilon,m}$ the matrix $\begin{psmallmatrix}
I_{\lfloor (m+1)/2\rfloor} & 0\\
0 & \epsilon I_{\lfloor m/2\rfloor}
\end{psmallmatrix}$.
\item[-] For any ring $R$ and $x\in \mathrm{\mathrm{Spec}}\ R$, we denote by $\kappa(x)$ the residue field $R_x/xR_x$ and by $\overline{\kappa(x)}$ an algebraic closure. 
\end{enumerate}

Now let us fix the integers $r$, $n$ with $r\leq \lfloor n/2 \rfloor$ and let $\e=\pm 1$. In this section, we introduce the matrix models $\mathcal{M}(n,r)$ and $\mathcal{M}(\e, n,r)$ for type A on the one hand, and types B, C, D on the other. They consist in giving $\mathbb{Z}$-group schemes $\G$, $\B$, $\T$, $\P$, $\L$, $\Z$, $\H$, morphisms $\Theta: \L\to \L$ and $\Pi:\Z\to \H$, elements $\mathfrak{e}$ and $\sigma_0$, and sets $\W$, $\W_{\P}$, $\W^{\P}$ as follows.

\subsection{Type A}

We define the matrix model $\mathcal{M}(n,r)$. For any ring $R$, we have on $R$-points:
\begin{align*}
\G(R)&=\G_n(R):=Gl_n(R),\\
\T(R)&=\T_n(R):= \begin{pmatrix}
* & & 0\\
 & \ddots & \\
0 & & *
\end{pmatrix} \subset \G(R),\\
\B(R)&=\B_n(R):= \begin{pmatrix}
* & * & *\\
 & \ddots & * \\
0 & & *
\end{pmatrix} \subset  \G(R),\\
\P(R)&=\P_{n,r}(R):=\Set{ \begin{pmatrix}
A & * & *\\
0 & B &* \\
0 & 0 & C
\end{pmatrix} \in \G(R)| \begin{smallmatrix}A,C\in \G_r(R)\\ B\in \G_{n-2r}(R)\end{smallmatrix}},\\
\L(R)&=\L_{n,r}(R):= \Set{ \begin{pmatrix}
A & 0 & 0\\
0 & B & 0 \\
0 & 0 & C
\end{pmatrix} \in \P(R) }= \Set{ \begin{pmatrix}
A & 0 & 0\\
0 & B & 0 \\
0 & 0 & C
\end{pmatrix}| \begin{smallmatrix}A,C\in \G_r(R)\\ B\in \G_{n-2r}(R)\end{smallmatrix}},\\
\Z(R)&=\Z_{n,r}(R):=\Set{ \begin{pmatrix}
A & * & *\\
0 & B & * \\
0 & 0 & A
\end{pmatrix} \in \P(R)},\\
\H(R)&=\H_{n,r}(R):=\Set{\begin{pmatrix}
A & 0 & 0\\
0 & 1 & 0 \\
0 & 0 & A 
\end{pmatrix}\in \P(R)},
\end{align*}
and 
\[\begin{array}{ccccc}
\Theta(R)=\Theta_{n,r}(R): & & \L(R)& \to & \L(R) \\
 & & \begin{pmatrix}
A & 0 & 0\\
0 & B & 0 \\
 0 & 0 & C
\end{pmatrix} & \mapsto & \begin{pmatrix}
C & 0 & 0\\
0 & B & 0 \\
 0 &0 & A 
\end{pmatrix},
\end{array}\]\\
\[\begin{array}{ccccc}
\Pi(R)=\Pi_{n,r}(R): & & \Z(R) & \to & \H(R) \\
 & & \begin{pmatrix}
A & * & *\\
0 & B &* \\
0 & 0 & A 
\end{pmatrix} & \mapsto & \begin{pmatrix}
A & 0 & 0\\
0 & 1 & 0 \\
0 & 0 & A 
\end{pmatrix}.
\end{array}\]
We also define
$$\mathfrak{e}=\mathfrak{e}_{n,r}:=\begin{pmatrix}
0 &0 & I_r\\
0& 0& 0\\
0& 0 &0
\end{pmatrix}$$
which is a square matrix of size $n$,
$$\s_{0}=\s_{0,n}:=\begin{pmatrix}
0 & & 1\\
& \iddots & \\
1 & & 0
\end{pmatrix}$$ which is a permutation in $\mathfrak{S}_n$, and finally we set
\begin{align*}
\W&=\W_n:=\mathfrak{S}_n,\\
\W_{\P}&:=\Set{\begin{pmatrix}
\sigma_1 & 0 & 0\\
0 & \s_2 & 0\\
 0 & 0& \s_3
\end{pmatrix}| \sigma_1,\ \s_3\in \mathfrak{S}_r,\ \s_2\in \mathfrak{S}_{n-2r}},\\
\W^{\P}&:=\Set{u\in \W | \begin{matrix}
u \mbox{ is increasing on }\\
\set{1,...,r},\ \set{r+1,...,n-r},\ \set{n-r+1,...,n}\end{matrix}}.
\end{align*} 

\subsection{Types B, C, D}\label{sec:matrix model BCD}

We define the matrix model $\mathcal{M}(\e,n,r)$ for integers $\e=\pm 1$ and $n\leq 2r$. Recall it encompasses the types $B$, $C$, $D$ for, respectively, $\e=1$ and odd $n$, $\e=-1$ and even $n$, $\e=1$ and even $n$. For $\G$ we adopt the picture presented in \cite{ChaputRomagny} in order to get smooth matrix groups over $\mathbb{Z}$. We explain some choices in the comments below. We need to introduce, for $m\in \mathbb{N}$, the orthogonal group $O_{2m}$ over $\mathbb{Z}$. We describe it as the closed $\mathbb{Z}$-subscheme of the affine space of $2m\times 2m$ matrices with equations
$${}^{\d}X
X=I_{2m} \quad \mbox{and} \quad  \sum_{k=1}^{ m }X_{k,i}X_{2m-k+1,i}=0,\ i\in \set{1,...,2m}$$ meaning any fiber of $O_{2m}$ is the group of linear transformations preserving the quadratic form \begin{equation}\label{eq:forme quadratic}
\sum_{k=1}^{ m }Y_{k}Y_{2m-k+1}.
\end{equation} Consider the determinant $det_{2m}$ on $O_{2m}$ as a regular function. There exists a unique regular function $\vartriangle_{2m}$, called the Dickson invariant, which satisfies $det_{2m}=1+2\vartriangle_{2m}$; in other words, $\vartriangle_{2m}$ takes values $0$ and $-1$ on the different components of $O_{2m}$ (see \cite[Lemma 4.1.4]{ChaputRomagny}). Let also $\dagger_{2m}$ be the difference between the $m$th and the $(m+1)$th vector of the usual basis of the $\mathbb{Z}$-free module $\mathbb{Z}^m$. We can now describe the data of our matrix model. For any ring $R$, we have on $R$-points: 

\begin{align*}
\G(R)&=\G_{\e,n}(R):=\left\lbrace\ \  \begin{matrix*}[l]
\Set{A\in O_n(R)|\vartriangle_{n}(A)=0} & \mbox{ if }\e=1\mbox{ and }n\mbox{ is even}\\
\\
 \Set{A\in O_{n+1}(R)| \begin{matrix} \vartriangle_{n+1}(A)=0\\
A\dagger_{n+1}=\dagger_{n+1}\end{matrix}}& \mbox{ if }\e=1\mbox{ and }n\mbox{ is odd}\\
 \\ 
 \Set{A\in M_n(R)| I_{-1,n} {}^{\d}A
I_{-1,n}A=I_n} & \mbox{ if }\e=-1\mbox{ and }n\mbox{ is even,} \end{matrix*}\right.\\
\T(R)&=\T_{\e,n}(R):=\left\lbrace\ \  \begin{matrix*}[l] \Set{\begin{psmallmatrix}
t_1 & & & & &\\
& \ddots & & & &\\
& & t_{n/2} & & &\\
& & & t_{n/2}^{-1} & &\\
& & & &\ddots & \\
& & & & & t_1^{-1}
\end{psmallmatrix}| \begin{smallmatrix} t_i\in R^*\\ \forall i\in\set{1,...,n/2}\end{smallmatrix}} & \mbox{ if }n\mbox{ is even} \\
\Set{\begin{psmallmatrix}
t_1 & & & & & & &\\
& \ddots & & & & & &\\
& & t_{(n-1)/2} & & & & &\\
& & & 1 & & & & \\
& & & & 1 & & &\\
& & & & & t_{(n-1)/2}^{-1} & &\\
& & & & & & \ddots & \\
& & & & & & & t_1^{-1}
\end{psmallmatrix}| \begin{smallmatrix} t_i\in R^*\\ \forall i\in \set{1,...,(n-1)/2}\end{smallmatrix}} & \mbox{ if }n\mbox{ is odd,}\end{matrix*}\right.\\
\B(R)&=\B_{\e,n}(R):= \Set{ \begin{pmatrix}
* & * & *\\
 & \ddots & *\\
0 & & *
\end{pmatrix} \in \G(R)},\\
\P(R)&=\P_{\e,n,r}(R):=\Set{ \begin{pmatrix}
A & * & *\\
0 & B & * \\
0 & 0 & C
\end{pmatrix} \in \G(R)| \begin{smallmatrix}A, C\in \G_{r}(R)\\B\in \G_{\e,n-2r}(R)\end{smallmatrix}}=\Set{ \begin{pmatrix}
A & * & *\\
0 & B & * \\
0 & 0 & {}^{\d}A^{-1}
\end{pmatrix} \in \G(R)| \begin{smallmatrix}A\in \G_{r}(R)\\B\in \G_{\e,n-2r}(R)\end{smallmatrix} },\\
\L(R)&=\L_{\e,n,r}(R):=\Set{ \begin{pmatrix}
A & 0 & 0\\
0 & B & 0 \\
0 & 0 & C
\end{pmatrix} \in\P(R) }= \Set{ \begin{pmatrix}
A & 0 & 0\\
0 & B & 0 \\
0 & 0 & ^{\delta}A^{-1}
\end{pmatrix} | \begin{smallmatrix}A\in \G_{r}(R)\\B\in \G_{\e,n-2r}(R)\end{smallmatrix} },\\
\Z(R)&=\Z_{\e,n,r}(R):=\Set{ \begin{pmatrix}
A & * & *\\
0 & B & * \\
0 & 0 & I_{-\e,r}AI_{-\e,r}
\end{pmatrix} \in \P(R)}=\Set{ \begin{pmatrix}
A & * & *\\
0 & B & * \\
0 & 0 & I_{-\e,r}AI_{-\e,r}
\end{pmatrix} \in \P(R)| \begin{smallmatrix}I_{-\e,r} {}^{\d}A
I_{-\e,r}A=I_r\\ B\in \G_{\e,n-2r}(R)\end{smallmatrix}},\\
\H(R)&=\H_{\e,n,r}(R):=\Set{\begin{pmatrix}
A & 0 & 0\\
0 & 1 & 0 \\
0 & 0 & I_{-\e,r} A I_{-\e,r} 
\end{pmatrix}  \in\P(R) | I_{-\e,r} {}^{\d}A
I_{-\e,r}A=I_r},
\end{align*}
and 
\[\begin{array}{ccccc}
\Theta(R)=\Theta_{\e,n,r}(R): & & \L(R)& \to & \L(R) \\
 & & \begin{pmatrix}
A & 0 & 0\\
0 & B & 0 \\
 0 & 0 & C
\end{pmatrix} & \mapsto & \begin{pmatrix}
I_{-\e,r}CI_{-\e,r} & 0 & 0\\
0 & B & 0 \\
 0 &0 & I_{-\e,r}AI_{-\e,r} 
\end{pmatrix},
\end{array}\]\\
\[\begin{array}{ccccc}
\Pi(R)=\Pi_{\e,n,r}(R): & & \Z(R) & \to & \H(R) \\
 & & \begin{pmatrix}
A & * & *\\
0 & B &* \\
0 & 0 & I_{-\e,r}A I_{-\e,r}
\end{pmatrix} & \mapsto & \begin{pmatrix}
A & 0 & 0\\
0 & 1 & 0 \\
0 & 0 & I_{-\e,r}A I_{-\e,r}
\end{pmatrix}.
\end{array}\]
We also define
$$\mathfrak{e}=\mathfrak{e}_{\e,n,r}:=\begin{pmatrix}
0 &0 & I_{-\e,r}\\
0& 0& 0\\
0& 0 &0
\end{pmatrix}$$
which is a square matrix of size $n$ for even $n$, of size $n+1$ for odd $n$,
$$ \s_{0}=\s_{0,\e,n}:=\left\lbrace\ \ \begin{matrix*}[l] \begin{psmallmatrix}
 0 & & & & & 1\\
 & & & & \iddots &\\
 & & 1 & 0 & &\\
 & & 0 & 1 & &\\
 & \iddots & & & &\\ 
 1 & & & & & 0
 \end{psmallmatrix} & \mbox{ if } \e=1,\ n\mbox{ is even and }n/2 \mbox{ is odd}\\
\begin{pmatrix}
 0 & &1\\
 & \iddots & \\
 1 & & 0
 \end{pmatrix} & \mbox{else}\end{matrix*}\right.$$
which is a permutation in $\mathfrak{S}_n$, and finally we set
\begin{align*}
\W&=\W_{\e,n}:=\Set{\s \in \mathfrak{S}_n | \begin{matrix} \check{\sigma}=\sigma \mbox{ and, if } n \mbox{ is even and } \e=1\mbox{, then }\\
\#\set{1\leq i \leq n/2 | \sigma (i)> n/2 } \mbox{ is even}
 \end{matrix}},\\
\W_{\P}&:=\Set{\begin{pmatrix}
\sigma & 0 & 0\\
0& v &0 \\
0&0 & \check{\sigma} \end{pmatrix} | \begin{matrix}
\sigma \in \mathfrak{S}_{r},\ v \in \mathfrak{S}_{n-2r},\\
\check{v}=v\mbox{ and, if } n \mbox{ is even and } \e=1\mbox{, then }\\
\#\set{1\leq i \leq (n-2r)/2 | v (i)> (n-2r)/2 } \mbox{ is even} \end{matrix}},\\
\W^{\P}&:=\Set{u\in \W |\begin{matrix}
u \mbox{ is increasing on }\set{1,...,r}, \set{r+1,...,\lfloor n/2 \rfloor}\\
 \mbox{ and } {u(\lfloor n/2 \rfloor)}<{u(\lfloor n/2 \rfloor+1+\frac{\e+(-1)^n}{2})}
\end{matrix}}.
\end{align*}

\subsection{Comments}\label{sec:comments matrix models}

Let us make some observations about the previous data. We prove the assumptions on flatness and smoothness below.  The assumptions concerning the reductiveness, semi-simpleness, the types and the Killing pairs\footnote{Note that at the level of schemes, all these notions demand smoothness, according to our setting which follows \cite[]{SGA3XIXaXXVI}.} then follow from the isomorphism Theorem (\cite[Theorem 1.1 Exposé XXV]{SGA3XIXaXXVI}) and the study of the geometric fibers of the involved groups. 

\paragraph{General}

\begin{enumerate}
\item[-] The couple $(\T,\B)$ is a Killing pair for $\G$.
\item[-] The nilpotent matrix $\mathfrak{e}$ can be identified with a section of the Lie algebra of $\G$ over $\mathbb{Z}$ and with a closed element of its Lie algebra over any geometric fiber. It is of order $2$ and rank $r$.
\item[-] We have $\Z=Z_{\G}(\mathfrak{e})$ the centralizer of $\mathfrak{e}$ in $\G$ for the adjoint action. Besides $\L\cap \Z=\L^{{\Theta}}$ and $\Z=\L^{\Theta}U_{\P}$ where $U_{\P}$ denotes the unipotent radical of $\P$. 
\item[-] The couple $(\T\times_{\G} \H,\B\times_{\G} \H)$ is a Killing pair for $\H$.
\item[-] The subgroup $\P$ of $\G$ is parabolic and $\L$ is its Levi subgroup containing $\T$.
\item[-] The Weyl group of $\G$ (respectively $\P$) related to $\T$ is naturally identified with $\W$ (respectively $\W_{\P}$). The set of minimum-length representatives of the quotient $\W/\W_{\P}$ then identifies with $\W^{\P}$. Besides, $\s_0$ is the longest element in $\W$.
\item[-] The isomorphism $\Theta$ is an involution of $\L$ which stabilizes $\T$ and $\B$ and $\Pi$ is a retraction from $\Z$ to $\H$.  

\end{enumerate}

\paragraph{Type A}
\begin{enumerate}
\item[-] The group $\G$ is reductive (not semi-simple) of type A$_{n}$.
\item[-] The group $\H$ is reductive (not semi-simple) of type A$_r$.
\item[-] On $W_{\P}$, $\Theta$ induces the bijection $$\begin{pmatrix}
\sigma & & \\
& v & \\
& & \sigma'
\end{pmatrix} \mapsto \begin{pmatrix}
\sigma' & & \\
& v & \\
& & \sigma
\end{pmatrix}.$$
\end{enumerate}

\paragraph{Types B, C, D}

\begin{enumerate}
\item[-] The group $\G$ is semi-simple of type B$_n$, D$_n$ in the cases $\e=1$ and odd $n$, $\e=1$ and even $n$ respectively. The base change $\G_{\mathbb{Z}_{(2)}}$ is semi-simple of type C$_n$ in the case $\e=-1$ and even $n$. In any case, the geometric fiber $\G_{\overline{x}}$ over any point $x\neq (2)\in \mathrm{\mathrm{Spec}}\ \mathbb{Z}$ is isomorphic to the group whose closed points consist in $\Set{A\in Sl_n(\overline{\kappa(x)})| I_{\e,n} {}^{\d}A
I_{\e,n}=A^{-1}}$. 
\item[-] We deduce that in the case $\e=1$, 
for any algebraically closed field $K$, the base change $\G_K$ is also isomorphic to the special orthogonal group presented in the introduction. Then $\mathfrak{e}$ identifies with a $n\times n$ matrix whose image is totally isotropic. 
\item[-] In the case of $\e=1$, $r$ is necessarily even. In the case of $\e=-1$, if $r$ is odd then we have an isomorphism $Z^0\rtimes \set{\pm 1} \simeq Z$ for any geometric fiber $Z:=\Z_{\overline{x}}$ over $x$.
\item[-] In the case of $\e=1$, $\H$ is semi-simple of type C$_r$. In the case of $\e=-1$, it is not even flat and its geometric fibers are not connected, but the base change $\H_{\mathbb{Z}_{(2)}}$ is smooth and the neutral components of its geometric fibers are semi-simple of type B$_r$.
\item[-] On $W_{\P}$, $\Theta$ induces the bijection $$\begin{pmatrix}
\sigma & & \\
& v & \\
& & \check{\s}
\end{pmatrix} \mapsto \begin{pmatrix}
\check{\s} & & \\
& v & \\
& & \sigma
\end{pmatrix}.$$
\end{enumerate}

We outline the following propositions.

\begin{prop}\label{prop:connexity}
In the models $\mathcal{M}(n,r)$ and $\mathcal{M}(1,n,r)$ (types A, B and D), $\Z$ and $\L^{\Theta}$ have geometric connected fibers. In $\mathcal{M}(-1,n,r)$ (type C) the fibers over $s\neq (2)\in \mathrm{\mathrm{Spec}}\ \mathbb{Z}$ have two connected components. However, the unipotent part $(\L^{\Theta}_{\overline{s}})_u$ is contained in the neutral component $(\L^{\Theta}_{\overline{s}})^0$.
\begin{proof}
It is clear for types A, B, D. For type C, and any $s\neq (2)\in \mathrm{\mathrm{Spec}}\ \mathbb{Z}$, we remark that $$\L^{\Theta}_{\overline{s}}\simeq O_{r{\overline{s}}}\times Sp_{n-2r\overline{s}}$$ and $$(\L^{\Theta}_{\overline{s}})^0\simeq SO_{r{\overline{s}}}\times Sp_{n-2r\overline{s}}.$$ This implies $$\L^{\Theta}_{\overline{s}} =(\L^{\Theta}_{\overline{s}})^0\cup \gamma(\L^{\Theta}_{\overline{s}})^0,$$ for a suitable $\gamma\in  \L^{\Theta}_{\overline{s}}$ with $det\ \gamma=-1$. Similar reasoning gives the same result for $\Z_{\overline{s}}$. We thus recognize two connected components. Since any unipotent matrix has its determinant equal to one, the above isomorphisms also yield the desired inclusion. 
\end{proof}
\end{prop}

\begin{prop}\label{prop:smoothness ZcapwB}
In the models $\mathcal{M}(n,r)$ and $\mathcal{M}(1,n,r)$ (types A, B, D), $\Z$, $\L^{\Theta}$ and $\Z\times_{\G}{}^w\B$ for any $w\in \W$ are smooth. In $\mathcal{M}(-1,n,r)$ (type C), they are not even flat but their base changes over $\mathbb{Z}_{(2)}$ are all smooth.  

\begin{proof}
We establish the claimed smoothness and flatness in the following way. They follow from the same general argument which also justifies the quite complicated picture of our groups in types B and D. It is based on the following lemma.

\begin{lem}\label{lem:smoothness criterium group scheme}
Let $S$ be a locally noetherian irreducible scheme and $G$ be a $S$-group scheme of finite type. Let $\eta$ denote the generic point of $S$. We assume that $G_{\overline{\eta}}$ is smooth. Then, for any $s\in S$ such that $\dim Lie(G_{\overline{s}})\leq \dim Lie(G_{\overline{\eta}})$, $G_{\overline{s}}$ is smooth. 

As a consequence, if the dimensions of the Lie algebras of the geometric fibers are the same and if $G$ is flat over $S$, then $G$ is smooth over $S$. 
\begin{proof}[Proof of Lemma \ref{lem:smoothness criterium group scheme}]
Replacing $G$ by its neutral component, we can assume that $G$ is irreducible. Let $s\in S$. Applying \cite[Lemma 13.1.1]{EGAIV8a15} to the dominant finite type morphism of irreducible schemes $G\to S$, we have the inequality of fiber dimensions 
$$\dim\ G_{\overline{\eta}}\leq \dim\ G_{\overline{s}}.$$ If we assume the smoothness of the geometric generic fiber and the inequality between the dimensions of Lie algebras, we thus have $$\dim\ Lie(G_{\overline{s}})\leq \dim\ Lie(G_{\overline{\eta}}) =\dim\ G_{\overline{\eta}}\leq \dim\ G_{\overline{s}}\leq \dim\ Lie(G_{\overline{s}})$$ which causes the smoothness of $G_{\overline{s}}$. The smoothness of $G$ over $S$ then follows from the smoothness fiber criterion for flat finite presentation morphisms. 
\end{proof}
\end{lem}

Let us now state our general argument. Let us consider the equations which define the desired group as a closed subscheme of the affine space of matrix $n\times n$. They are all linear or quadratic, with the exception of the one concerning the vanishing of the Dickson invariant. Some of the quadratic ones are given by 
\begin{equation}\label{eq:torsion}
X{}^{\delta}X=I_{2m}
\end{equation}
for some integers $m$ and indeterminates $2m\times 2m$ matrices $X$. They lead to $2$-torsion with the equations depending on $i\in \set{1,...,2m}$  $$\sum_{k=1}^{ 2m }X_{ki}X_{2m-k+1,i}=\delta_{i,2m-i+1}$$ namely $$2\sum_{k=1}^{ m }X_{ki}X_{2m-k+1,i}=\delta_{i,2m-i+1}.$$ That can be compensated by refining them with \begin{equation}\label{eq:torsion compensated}
\sum_{k=1}^{ m }X_{ki}X_{2m-k+1,i}=0.
\end{equation} This problem does not occur for the other quadratic equations arising from \begin{equation}\label{eq:non torsion}
I_{-1,m}{}^{\delta}XI_{-1,m}X=I_m.
\end{equation} Note that the Dickson invariant does not imply torsion whereas its existence yields one for the equation $det_{2m}=1$. Hence, the group at stake is flat or non-flat over $\mathbb{Z}$ according to whether (\ref{eq:torsion}) or (\ref{eq:non torsion}) is involved and potentially compensated. In any case, it is flat over $\mathbb{Z}_{(2)}$. On the other hand, it is finitely presented, and we can check that the Lie algebra of any of its geometric fibers is of constant dimension. Since the generic geometric fiber is smooth as an algebraic group in characteristic zero, we deduce the smoothness of all geometric fibers by the Lemma \ref{lem:smoothness criterium group scheme}. Let us emphasize that the constancy of the Lie algebra dimension follows from the parity of $2m$. Indeed, it prevents the presence of squared terms in (\ref{eq:torsion}) and (\ref{eq:torsion compensated}), so that the differentials do not bring $2$-torsion and therefore a leap of dimension. This explains why we describe $\G_{1,n}$ into $\G_{1,n+1}$ for odd $n$. 
\end{proof}
\end{prop}

\section{Birationality in types A, B, C and D}\label{sec:birationality}

We start with proving Theorem \ref{theo:1 birational rational} except for the claim on rationality. Our conclusions include a partial version of the result for type C; see Section \ref{sec:conclusion birationality}.

\subsection{Two lemmas}

We will use the two basic and technical lemmas below.

\begin{lem}\label{lem:équivalence birationnelle cas général}
Let $G$ be a connected reductive algebraic group over $\k$ and $(T,B)$ be a Killing pair. Let $H\subset Z\subset G$ be connected subgroups equipped with a retraction $\varpi: Z\to H$ and such that $B_H:=B\cap H$ is a Borel subgroup of $H$.  Let $w$ be in $W$, the Weyl group of $G$. We assume that:
\begin{enumerate}
\item $\dim Z/\left(Z\cap \conj{w}{B}\right) =\ell(w)+\dim H/B_{H}$
\item $z^{-1}\varpi (z) \in \overline{\conj{w}{B}B}$ for each $z\in Z$
\item $Z\cap \conj{w}{B}\subset \varpi^{-1}(B_H) (Z\cap{}^wB)^0$
\item The scheme-theoretic intersection $Z\times_G\conj{w}{B}$ is reduced.
\end{enumerate} Then $Z$ acts on $H\times^{B_H}\overline{B\cdot wB}$ by $z\cdot [h,gB]=[\varpi(z)h,h^{-1}\varpi(z)^{-1}zhgB]$ and the map
\begin{equation}
\pi: \begin{array}{cccc}
H\times^{B_H} \overline{B\cdot wB} & \to & \overline{Z\cdot wB}\\

[h,gB] & \mapsto & hgB.
\end{array}
\end{equation} is birational and $Z$-equivariant. 
\end{lem}

\begin{lem}\label{lem:formule dimension sous-groupe symétrique} Let $G$, $B$, $T$, $W$, $H$, $Z$, $B_H$, $\varpi$ be as in the previous lemma. Let $P\supset B$ be a parabolic subgroup with $U_P$ its unipotent radical and $L$ its Levi subgroup containing $T$. Let $\theta$ be an involution of $L$ which stabilizes $B\cap L$ and $T$. We assume that $Z$ is the subgroup $(L^{\theta})^0U_P\subset P$ and we fix $w\in W$. 
\begin{enumerate}
\item[(i)] If the subvariety of unipotents $L^{\theta}_u$ of $L^{\theta}$ is contained in $L_Z$, then $$Z\cap \conj{w}{B} \subset  \varpi^{-1}(B_H)(Z\cap \conj{w}{B})^0.$$
\item[(ii)] If $char(\k)\neq 2$ and if $w=\tau v$ is the decomposition of $w$ in $W_P(W^P)^{-1}$, then
$$  \dim Z/Z\cap \conj{w}{B} = \ell(w)+\dim Z/Z\cap B +\ell(\tau ^{-1}\theta (\tau))/2-\ell(\tau).$$
\end{enumerate}
\end{lem}

\subsubsection{Proof of Lemma \ref{lem:équivalence birationnelle cas général}}

Let $X$ be $\overline{Z\cdot wB}$ and $\hat{X}$ be the subvariety $(\iota\times id)^{-1}(\overline{G\cdot (eB,wB)})$ of $H/B_H\times G/B$ where $\iota$ is the immersion $H/B_{H} \hookrightarrow G/B$.  We have the isomorphism \begin{equation}\label{eq:iso hatX produit contracté}
\hat{X}\simeq H\times^{B_H} \overline{B\cdot wB}
\end{equation} over $H/B_H$ as $H$-equivariant bundles and over $G/B$ thanks to $\pi$ and the second projection $pr_2:H/B_H\times G/B\to G/B$. Considering this isomorphism, it suffices to show that $pr_2$ induces a well defined birational morphism $\hat{X}\to X$ and, transporting the action, to ensure that $\hat{X}$ is stable for the $Z$ action on $H/B_H\times G/B$ given by $z\cdot (hB_H,gB)=(\varpi(z)hB_H,zgB)$. We will proceed in several steps.
\begin{enumerate}
\item With the hypothesis $1$, we have $\dim \hat{X}=\ell (w)+\dim H/B_H=\dim Z\cdot wB=\dim X$. 
\item If $f$ denotes the morphism $g\mapsto (gB,wB)$,  we have for all $b'\in {}^wB$ and $b\in B$, $f(b'b)=(b'bB,wB)=(b'B,b'wB)=b'\cdot (B,wB)$. Thus $f(^wBB)\subset \overline{G\cdot (eB,wB)}$ and $f(\overline{^wBB})\subset \overline{G\cdot (eB,wB)}$.  Thanks to the hypothesis $2$, we have then for all $z\in Z$ $$\iota \times id\ (z\cdot (eB_{H}, wB))=\iota \times id\ (\varpi(z)B_{H},zwB)=z\cdot f(z^{-1}\varpi(z))\in z\cdot f(\overline{^wBB})\subset \overline{G\cdot (eB,wB)}.$$ We deduce that $Z\cdot (eB_H,wB)\subset \hat{X}$.

\item The map $pr_2$ induces a surjective $Z$-equivariant morphism $Z\cdot (eB_H,wB) \to Z\cdot wB$. The previous points give then $\dim Z\cdot wB \leq \dim Z\cdot (eB_H,wB) \leq \dim \hat{X} =\dim Z\cdot wB$ so $\dim Z\cdot wB =\dim Z\cdot (eB_H,wB)$. This also implies $\dim Z\cap {}^wB=\dim Z\cap {}^wB\cap \varpi ^{-1}(B_{H})$ so that $\left( {}Z\cap {}^wB\right)^0=\left( {}Z\cap {}^wB\cap \varpi ^{-1}(B_{H})\right)^0$ and $(Z\cap {}^wB)^{0}\subset \varpi ^{-1}(B_{H})$.  With hypothesis $3$, we deduce $\varpi(Z\cap {}^wB)\subset B_{H}$.  Therefore 
\begin{align*}
pr_2^{-1}(wB)\cap Z\cdot (eB_H,wB)&=\Set{(\varpi(z)B_{H},zwB) |\ z\in Z, zwB=wB}\\
&=\Set{(\varpi(z)B_{H},zwB)| \ z\in Z\cap \conj{w}{B}}\\ &=\varpi(Z\cap {}^wB)\cdot B_{H}\times \lbrace wB\rbrace\\
&= \lbrace (eB_{H},wB)\rbrace,
\end{align*} and $Z\cdot (eB_H,wB)\to Z\cdot wB$ is bijective. But it is also separable thanks to hypothesis $4$. It is finally an isomorphism by the Zariski Main Theorem. 

\item The previous points imply in particular that $\dim \hat{X}=\dim Z\cdot (eB_H,wB)$ and then $\overline{Z\cdot (eB_H,wB)}=\hat{X}$. Hence, $Z$ preserves $\hat{X}$ as desired.  

\item To conclude, note that we have $Z\cdot wB\subset pr_2(\hat{X})$. Since $H/B_{H}$ is complete, $pr_2(\hat{X})$ is closed in $G/B$ and we deduce $X\subset pr_2(\hat{X})$. Irreducibility and dimension formula give then $\dim pr_2(\hat{X})\leq \dim \hat{X}=\dim X$. Therefore $\hat{X}\to X$ is well defined and surjective.
\end{enumerate}

We have all the desired statements, with an isomorphism between the dense open orbits $Z\cdot (eB_H,wB)$ and $Z\cdot wB$.

\subsubsection{Proof of Lemma \ref{lem:formule dimension sous-groupe symétrique}}

Let $\Phi_L$ be the set of roots of $L$ and put $L_Z:=(L^{\theta})^0$, $B_L:=B\cap L$ and $T_Z:=(T\cap Z)^0=(T\cap L_Z)^0=T^{\theta,0}$. The involution $\theta$ and the elements of $W_P$ act linearly on the vector space $\mathbb{R}\otimes _{\mathbb{Z}}\Phi_L$ and let us denote these respective actions by $\star$ and $\diamond$. Note that $B_L$ is a $\theta$-stable Borel subgroup of $L$. Besides, since $T\subset B_L$ is a $\theta$-stable maximal torus of $L$, $T_Z$ is a regular subtorus of $L$ (see for example \cite[Lemma 4]{BrionHelminck}\footnote{Which is also valid for characteristic two.}) and it follows that it is a maximal torus of $L_Z$ and of $Z$.

\begin{enumerate}
\item Let us first prove $(i)$. Since $T\subset L$, there exists a cocharacter $\lambda: \mathbb{G}_m\to T$ such that $U_P=\Set{ x \in G| lim_{a\to 0}\lambda (a)x\lambda (a)^{-1}=1 }$ and $L=Z_G(Im\ \lambda)$.  Hence, if $z=lv\in Z\cap \conj{w}{B}$ with $l\in L_Z$ and $v\in U_P$ then $\lambda (a)z\lambda (a)^{-1}\in \conj{w}{B}$ for all $a\in \mathbb{G}_m$ and $\lambda (a)z\lambda (a)^{-1}=l \lambda (a)v\lambda (a)^{-1}\to l$ when $a\to 0$. We therefore have $l,v\in \conj{w}{B}$ and $z\in (L_Z\cap \conj{w}{B})(\conj{w}{B}\cap U_P)$ so that \begin{equation}\label{eq:inclusion2 pour ZcapwB}
Z\cap \conj{w}{B} \subset (L_Z\cap \conj{w}{B})(Z\cap \conj{w}{B})_u.
\end{equation}
But we have
\begin{equation}\label{eq:inclusion0 pour ZcapwB}
L_Z\cap \conj{w}{B} \subset (T\cap Z)(Z\cap \conj{w}{B})_u.
\end{equation} To see this pick $x\in L_Z\cap \conj{w}{B}$ and let $x=tv$ be its decomposition with $t\in T$ and $v\in (^wB)_u$ in the connected solvable group ${}^wB$. We have $v=t^{-1}x\in TL_Z\subset L$ and we can apply $\theta$ on $v$. Since $x\in L_Z$ we have $tv=\theta (t)\theta (v)$ and $\theta (v)=\theta (t)^{-1}tv\in \conj{w}{B}$ because $\theta$ stabilizes $T\subset \conj{w}{B}$. Hence the unipotent element $\theta (v)$ is in the group $(\conj{w}{B})_u$ and $\theta (t)^{-1}t=\theta (v)v^{-1}$ is also unipotent. It is thus the unit element which implies $t,v\in L^{\theta}$. By hypothesis $L^{\theta}_u\subset L_Z$ so that $v$, and then $t$, are in $L_Z$. Therefore $x\in (T\cap L_Z)((\conj{w}{B})_u\cap L_Z)\subset (T\cap Z)(Z\cap \conj{w}{B})_u$.  We have also
\begin{equation}\label{eq:inclusion1 pour ZcapwB}
T\cap Z\subset \varpi^{-1}(B_H).
\end{equation} Indeed, since $\varpi:Z\to H$ is a retraction, we have $B_H=B\cap H\subset B\cap Z$ and $B_H=\varpi(B_H)\subset \varpi(B\cap Z)$ which is an equality because $\varpi(B\cap Z)$ is solvable. 
We have also 
\begin{equation}\label{eq:inclusion unipotent composante connexe}
(Z\cap \conj{w}{B})_u\subset (Z\cap \conj{w}{B})^0.
\end{equation} 
Indeed, the unipotent subgroup $(Z\cap \conj{w}{B})_u$ of the connected group $Z$ is contained in a Borel subgroup $\overset{\sim}{B}_Z$ of $Z$ (see for example \cite[Theorem 30.4]{Humphreys}). Since it is stable under the action of the maximal torus $T_Z$, we can assume that $T_Z\subset \overset{\sim}{B}_Z$ by the Borel fixed point theorem. As for $P$ and $U_P$, the unipotent radical $(\overset{\sim}{B}_Z)_u$ is described by $\Set{x\in Z| lim_{a\to 0}\mu (a)x\mu (a)^{-1}=1 }$ where $\mu$ is a suitable cocharacter $\mathbb{G}_m\to T_Z$. Therefore, conjugating by $\mu(a)$ and taking the limit when $a\to 0$, any $x\in (Z\cap \conj{w}{B})_u$ can be contracted, inside $Z\cap \conj{w}{B}$, to the unit element, and thus, belongs to $(Z\cap \conj{w}{B})^0$.

Combining (\ref{eq:inclusion2 pour ZcapwB}), (\ref{eq:inclusion0 pour ZcapwB}), (\ref{eq:inclusion1 pour ZcapwB}),  (\ref{eq:inclusion unipotent composante connexe}) we obtain the desired inclusion.

\item Let us now prove $(ii)$ under the assumption $char(\k)\neq 2$. On the one hand we can use a result of Richardson and Springer (\cite[Proposition 3.3.4]{RichardsonSpringer1993} reformulating \cite[Theorem 4.6]{RichardsonSpringer1990}) on the involution fixed-points subgroup $L^{\theta}$ of $L$ and we get (replacing harmlessly $L^{\theta}$ with $L_Z$):
$$\dim L_Z/L_Z\cap {}^{\tau}B=\dim L_Z/L_Z\cap B+\left(\ell\left(\tau^{-1}\theta(\tau)\right)+\dim E_-\left(\tau^{-1}\theta(\tau)\right)-\dim E_-(1)\right)/2$$
where $E_-(\sigma):=\Set{v\in \mathbb{R}\otimes _{\mathbb{Z}}\Phi_L \vert \ \sigma\diamond (\theta\star v)=-v }$ for any $\sigma\in W_P$. But for all $v\in \mathbb{R}\otimes _{\mathbb{Z}}\Phi_L$ we have $\theta(\tau)\diamond(\theta \star v)=\theta \star (\tau \diamond v)$ whence the equivalences
\begin{align*}
&(\tau^{-1}\theta(\tau))\diamond (\theta \star v)=-v\\
\Leftrightarrow &\ \tau\diamond \tau^{-1} \diamond \theta(\tau)\diamond (\theta \star v)=-\tau \diamond v \\
\Leftrightarrow &\ \theta(\tau)\diamond (\theta \star v)=-\tau \diamond v \\
\Leftrightarrow &\ \theta \star (\tau \diamond v)=-\tau \diamond v
\end{align*} which imply $E_-\left(\tau^{-1}\theta(\tau)\right)=\tau^{-1}\diamond E_-(1)$. Besides, since $U_P\subset B$, we have $L_Z/L_Z\cap B\simeq Z/Z\cap B$. We thus have 
\begin{equation}\label{eq:egalité dimension finale avec Richardson Springer}
\dim L_Z/L_Z\cap {}^{\tau}B=\dim Z/Z\cap B+\ell\left(\tau^{-1}\theta(\tau)\right)/2.
\end{equation}

In addition, since $char(\k)\neq 2$, the fixed point subgroup $L^{\theta}$ has a finite number of orbits in the flag variety $L/B_L$ (see \cite[§4]{Springer1985}). Thus $L_Z$ is spherical in $L$. We can then apply a result\footnote{The article is written for characteristic zero but we only need here the second part, which is valid for any characteristic.} of Chaput, Fresse and Gobet (\cite[Theorem 7.2 (c)]{ChaputFresse}) for the subgroup $Z$ constructed by parabolic induction $Z=L_ZU_P$. We get, for $\tau\in W_P$ and $\upsilon \in W^P$:
$$\dim Z/Z\cap {}^{\tau  \upsilon}B=\ell( \upsilon)+ \dim L_Z/L_Z\cap {}^{\tau}B.$$
But with the decomposition $w=\tau \upsilon$ we have $\ell(\upsilon)=\ell(w)-\ell(\tau)$ and therefore 
\begin{equation}\label{eq:dimension finale avec Chaput Fresse Gobet}
\dim Z/Z\cap {}^{w}B=\ell(w)-\ell(\tau) + \dim L_Z/L_Z\cap {}^{\tau}B.
\end{equation}
Combining (\ref{eq:egalité dimension finale avec Richardson Springer}) and (\ref{eq:dimension finale avec Chaput Fresse Gobet}) we obtain the desired formula.
\end{enumerate}

\subsection{Applying to matrix models}\label{sec:application through matrix models}

Let us now fix the framework presented in the introduction for $G$, $e$ and $Z$ relative to types A, B, C, D. Recall that we assume that the image of $e$ is totally isotropic for types B and D (which is automatically satisfied if $char\ k\neq 2$) and that we only consider $char(\k)\neq 2$ for type C. Let $r$ be the rank of $e$. From the comments in Section \ref{sec:matrix models}, we can therefore assume that $G$ is the base change $\G_{\mathds{k}}$ of $\G$ described in the matrix model $\mathcal{M}(n,r)$ or $\mathcal{M}(\e,n,r)$.

In addition, we know for ${Gl_n}_{\k}$ that Young diagrams parametrize nilpotent orbits and that those related to elements of order two are characterized solely by the rank of these elements (two-column case). For $char(\k)\neq 2$ this is also the case for the ${O_n}_{\k}$ and ${Sp_n}_{\k}$ groups (see, for example, \cite[§1.6]{Jantzen}). For $char(\k)=2$ and types B, D, if we want to characterize the orbit of a nilpotent element $N$, then we must add to the rank the datum of the sequence $(\chi_m^N)_{m\geq 1}$ whose $m$th term is $$\chi_m^N = min\set{l\geq 0| N^l\left(KerN^m\right)\mbox{ is totally isotropic}}$$(see \cite[Theorem 3.8]{Hesselink1979}). But such a sequence is totally determined by the rank $s$ of $N$, as soon as $N$ is of order two with a totally isotropic image. Indeed, in this situation, we first obviously have $\chi_1^N\leq 1$ and $\chi_m^N=1$ for all $m\geq 2$. Moreover, we have $2s\leq n$ since $ImN\subset KerN$. If $2s=n$, the rank theorem implies that $KerN=ImN$ so that $KerN$ is totally isotropic and $\chi_1^N=0$. If $2s<n$, the dimension of $Ker N$ exceeds $n/2$ so that it cannot be a totally isotropic subspace and we deduce $\chi_1^N\geq 1$. 

Thanks to an inner conjugation (for types A, C) or (potentially) an outer conjugation (for types B, D), we can therefore assume that $e$ and then $Z$ respectively coincide with $\mathfrak{e}$ and the base change $\Z_{\mathds{k}}$ of the appropriate matrix model. We supplement this data with $B$, $T$, $P$, $L$, $H$, $B_{H}$, $T_{H}$, $\theta$, $\varpi$ and $W$, $W_P$, $W^P$, which stand for the various base changes $\B_{\mathds{k}}$, $\T_{\mathds{k}}$, $\P_{\mathds{k}}$, $\L_{\mathds{k}}$, $(\H_{\mathds{k}})^0$, ${(\B\times_{\G}\H)}_{\mathds{k}}$, ${(\T\times_{\G}\H)}_{\mathds{k}}$, $\Theta_{\mathds{k}}$, $\Pi_{\mathds{k}}$ and the sets $\W$, $\W_{\P}$, $\W^{\P}$. We also set $L_Z:=(L^{\theta})^0$, $B_L:=B\cap L$. We keep the notation $\sigma_0$ for the longest element in $W$ and denote by $u_0$ the permutation $\sigma_{0,n-2r}$ (type A), or $\sigma_{0,\e,n-2r}$ (types B, C, D).

For any integer $m$, and $\epsilon=\pm 1$ we will also denote more generally by $G_m$, $B_m$, $T_m$, $G_{\epsilon, m}$, $B_{\epsilon, m}$, $T_{\epsilon, m}$, $W_m$, $W_{\e,m}$ the base changes ${\G_{m}}_{\mathds{k}}$, ${\B_{m}}_{\mathds{k}}$, ${\T_{m}}_{\mathds{k}}$, ${\G_{\epsilon, m}}_{\mathds{k}}$, ${\B_{\epsilon, m}}_{\mathds{k}}$, ${\T_{\epsilon, m}}_{\mathds{k}}$, and the sets $\mathcal{W}_{m}$, $\mathcal{W}_{\e,m}$.  We also introduce the natural maps $$\zeta_{m}: N_{G_m}(T_m)\to W_m$$ which map a monomial matrix to the permutation (matrix) it induces.

In order to apply the previous lemmas to this situation, we begin by proving several properties for the type A and then for the three types B, C, D.

\paragraph{Restrictions for type C} In type C, we assume that $char(\k)\neq 2$, that $r$ is odd and we replace $Z$ and $H$ by their neutral components $Z^0$ and $H^0$ (see Remark \ref{rqe:type C}).  

\subsubsection{Type A}

We begin with a proposition of independent interest.

\begin{prop}\label{prop:type A T fixed point}

We have the following parametrization of the $Z$-orbits in $G/B$

$$\begin{array}{ccccl}
& & W_r\times W^P & \to & Z \backslash (G/B)\\ 
& & (u,v) & \mapsto & {Z}\cdot \begin{psmallmatrix} u & &\\ 
& 1 & \\ & & 1 \end{psmallmatrix}v^{-1}B\end{array}$$ 
As a consequence, the number of $Z$-orbits is finite. Besides, they all possess at least one $T$-fixed point. 
\begin{proof}
Let us recall that the Bruhat decomposition implies a bijection 
$$\begin{array}{ccccl}
& & W_r & \to & G_r\backslash \left(G_r/B_r\times G_r/B_r\right)\\ 
& & u  & \mapsto & G_r\cdot (\s B_r,eB_r).\end{array}$$ 
But, considering the block shapes of $L$ and $L_Z$, there is a bijection 
$$\begin{array}{ccccl}
& & G_r\backslash \left(G_r/B_r\times G_r/B_r\right) & \to & L_Z\backslash \left(L/B_L\right)\\ 
& & G_r\cdot (xB_r,yB_r) & \mapsto & L_{Z}\cdot \begin{psmallmatrix} x & &\\ & 1 & \\ & & y \end{psmallmatrix}B_{L}.\end{array}$$ 
We therefore have a bijection
$$\begin{array}{ccccl}
& & W_r & \to & L_Z\backslash \left(L/B_L\right)\\ 
& & u & \mapsto & L_{Z}\cdot \begin{psmallmatrix} u & &\\ & 1 & \\ & & 1 \end{psmallmatrix}B_{L}.\end{array}$$ 
We can then apply again \cite[Theorem 7.2 (a)]{ChaputFresse} to the spherical subgroup $Z=L_ZU_P$ which gives finally the desired parametrization $$\begin{array}{ccccl}
& & W_r\times W^P & \to & Z \backslash (G/B)\\ 
& & (u,v) & \mapsto & {Z}\cdot \begin{psmallmatrix} u & &\\ 
& 1 & \\ & & 1 \end{psmallmatrix}v^{-1}B.\end{array}$$
The assertion on $T$-fixed points follows easily. 
\end{proof}
\end{prop}

Now, we consider the three essential results below.
 
\begin{prop} \label{prop:type A corriger par z}
Let $v\in W$. Then, there exists $z_0\in Z\cap N_G(T)$ such that $w:=\zeta_n (z_0)^{-1}v$ satisfies
\begin{enumerate}
\item[$(a)$] $w^{-1}$ induces $u_0$ on $\set{r+1,...,n-r}$; in other words, $w^{-1}$ is decreasing on this set.
\item[$(b)$] $w^{-1}$ is increasing on $\set{1,...,r}$.
\end{enumerate}
\begin{proof}
If $v\in W$, there exists $\s_1\in \mathfrak{S}_{n-2r}$ and $\s_2\in \mathfrak{S}_r$ such that $v^{-1}\begin{pmatrix}
1 & &\\
& \s_1 &\\
& & 1
\end{pmatrix}$ decreases on $\set{r+1,...,n-r}$ and $v^{-1}\s_2$ increases on $\set{1,...,r}$. We fix then 
$$z_0:=\begin{pmatrix} \s_2 & & \\ & \s_1 & \\ & & \s_2 \end{pmatrix}$$
which is obviously an element of $Z\cap N_G(T)$, and which is also equal, as a matrix, to $\zeta_n(z_0)$ in $W$. Hence $v^{-1}\zeta_n(z_0)$ is $v^{-1}\s_1$ on $\set{r+1,...,n-r}$ and is $v^{-1}\s_2$ on $\set{1,...,r}$,  so that $\zeta_n (z_0)^{-1}v$ is the desired element.
\end{proof}
\end{prop}

 \begin{prop}\label{prop:type A a imply hypo 2 }
Suppose that $w\in W$ satisfies property $(a)$ of Proposition \ref{prop:type A corriger par z}. Then $z^{-1}\varpi (z)\in \overline{\conj{w}{B}B}$ for all $z\in Z$.
\begin{proof}
Let $w\in W$ satisfying the assumption.  
Since $w$ induces $u_0$ we get \begin{equation} \label{eq:1 condition i theo type A}
 \begin{Large} 
 \begin{psmallmatrix}
1 &  & \\
 & ^{u_0} B_{n-2r} & \\
 &  & 1
\end{psmallmatrix} \subset \ {}^wB.\end{Large} 
\end{equation}
But $u_0$ is the longest element in $W_{n-2r}$ so that $${}^{u_0}B_{n-2r}B_{n-2r}\simeq B_{n-2r}u_0B_{n-2r}$$ is dense in $G_{n-2r}$ and (\ref{eq:1 condition i theo type A}) thus implies
\begin{equation} \label{eq:2 condition i theo type A}
\begin{pmatrix}
1 & & \\
 & G_{n-2r} & \\
 & & 1
\end{pmatrix}\subset \overline{{}^wBB}.
\end{equation} Let now $z\in Z$. The element $z^{-1}\varpi (z)$ has the shape
$$\begin{pmatrix}
1 & A_1 & A_2\\
0 & C & A_3\\
0 & 0 &1 
\end{pmatrix}=\begin{pmatrix}
1 & & \\
 & C & \\
 & &1 
\end{pmatrix} \begin{pmatrix}
1 & A_1 & A_2\\
0 & 1 & C^{-1}A_3\\
0 & 0 &1 
\end{pmatrix}$$ with $C\in G_{n-2r}$ and suitable matrices $A_i$. We conclude that $z^{-1}\varpi (z)\in \overline{^wBB} B=\overline{^wBB}$ by (\ref{eq:2 condition i theo type A}).
\end{proof}
\end{prop}

\begin{prop}\label{prop:type A a b quasi imply hypo 1}
Suppose $w\in W$ satisfies properties $(a)$ and $(b)$ of Proposition \ref{prop:type A corriger par z}. Let $w=\tau\nu$ be the decomposition of $w$ in $W_P(W^P)^{-1}$. Then $$\ell (\tau) -(\ell (\tau ^{-1}\theta (\tau))/2=\dim Z/Z\cap B-\dim H/B_H.$$
\begin{proof}
Let $w\in W$ satisfying the assumptions. There exists $\s\in \mathfrak{S}_r$ such that $w^{-1}\begin{psmallmatrix}
1 & &\\
& 1 &\\
& & \s
\end{psmallmatrix}$ is increasing on $\set{n-r+1,...,n}$. We fix then
$$\tau:=\begin{psmallmatrix}
1 & &\\
& u_0 & \\
& & \s
\end{psmallmatrix}$$ which is in $W_P$. Besides, $w^{-1}\tau$ is in  $W^P$ because it restricts to $w^{-1}$ on $\set{1,...,r}$, to $w^{-1}\begin{psmallmatrix}
1 & &\\
& u_0^{-1} & \\
& & 1
\end{psmallmatrix}$ on $\set{r+1,...,n}$ and to $w^{-1}\begin{psmallmatrix}
1 & &\\
& 1 &\\
& & \s
\end{psmallmatrix}$ on $\set{n-r+1,...,n}$, which are all increasing permutations by assumptions on $w$. Besides, we have
$$ \tau ^{-1}\theta (\tau)=\begin{psmallmatrix}
1 & & \\
 & u_0^{-1} & \\
& & {\sigma}^{-1} 
\end{psmallmatrix} \begin{psmallmatrix}
{\sigma} & & \\
 & u_0 & \\
 & & 1
\end{psmallmatrix}
= \begin{psmallmatrix}
\sigma & & \\
 & 1 & \\
 & & {\sigma}^{-1} 
\end{psmallmatrix}.$$ We thus have $\ell (\tau ^{-1}\theta (\tau))=2\ell(\s)$. Besides, $\ell (\tau)=\ell(u_0)+\ell (\s)$. But an easy computation gives $\dim Z/Z\cap B-\dim H/B_H=\dim G_{n-2r}/B_{n-2r}=\ell (u_0)$. Combining all these equalities, we obtain the desired formula. 
\end{proof}
\end{prop}

\subsubsection{Types B, C and D}\label{sec:B,C,D}

The previous properties have their analogs in types B, C and D. Before giving them, let us begin with some material adapted to this setting.

\paragraph{Some preliminaries}

For $v\in W$, we consider the quantity $$d_v:=\#\set{r+1\leq i\leq \lfloor n/2 \rfloor | v^{-1}(i)>\lfloor n/2 \rfloor}.$$ Then define $s_{v}\in \mathfrak{S}_{n-2r}$ as $$s_v:=\left\lbrace \begin{matrix} \left( \frac{n-2r}{2}\ \ \frac{n-2r}{2}+1\right) & \mbox{if }\e=1,\ n \mbox{ is even, } d_v \mbox{ is odd }  \\
id & \mbox{else.}\end{matrix}\right. $$ We motivate our definition by the following fact: 

\begin{prop}\label{prop:fait sv}
If $u\in \mathfrak{S}_{n-2r}$ is induced by $v^{-1}$ on $  \set{r+1,...,n-r}$ then $s_{v}u \in W_{\e,n-2r}$. 
\begin{proof}
Indeed, let $u$ be such an element. Since $\widecheck{v^{-1}}=v^{-1}$, we have $\check{u}=u$ and thus $\widecheck{s_v u}=s_v u$ in $\mathfrak{S}_{n-2r}$, so that it is enough to show the parity of $q:=\#\set{1\leq i\leq (n-2r)/2| s_v u(i)>(n-2r)/2}$ for $\e=1$ and even $n$. In this case, $i\mapsto r+i$ identifies $\set{1\leq i\leq (n-2r)/2| u(i)>(n-2r)/2}$ with $\set{r+1\leq i\leq n/2 | v^{-1}(i)> n/2}$. Therefore, if $d_v$ is even then $s_v=id$ and $q=d_v$; if $d_v$ is odd, then $s_v=\left( \frac{n-2r}{2}\ \ \frac{n-2r}{2}+1\right)$ and $q=d_v\pm 1$. In each case $q$ is even and we have proved the proposition.
\end{proof}
\end{prop}

We will also use the following results concerning length and permutations.

\begin{prop}\label{prop:préliminaire technique types B,C,D}
Let $d\in \mathbb{N}$, $\epsilon=\pm 1$, $v\in \mathfrak{S}_d$ and $u\in W_{\epsilon,d}$. We have 
\begin{enumerate}
\item[(i)] $\ell_{d}(u)=2\ell(u)+\epsilon\#\set{1\leq i\leq \lfloor d/2 \rfloor | u (i)>\lfloor d/2 \rfloor}$.
\item[(ii)] If $v(i)<v (j)$ or $v (\overline{i})>v (\overline{j})$ for any $1\leq i<j\leq d$, then $\ell_d(\check{v}v ^{-1})=2\ell_d(v)$.
\item[(iii)] There exists $\sigma\in \mathfrak{S}_{d}$ such that $\check{\sigma}=\sigma$ and $v\sigma (i)<v\sigma (j)$ or $v\sigma (\overline{i})>v\sigma (\overline{j})$ for any $1\leq i<j\leq d$.
\end{enumerate}
\begin{proof}
We prove the assertions separately. 
\begin{enumerate}
\item[$(i)$] It is the result of taking proper account of the root systems associated with the different types.
\item[$(ii)$] We actually have 
\begin{equation}\label{eq:longueur sigmacheck sigmamoinsun}
\ell_d(\check{v}v^{-1})=2\ell_d(v)-2\#\set{1\leq i<j\leq d | v (i)>v (j) \text{ and } v (\overline{i})<v (\overline{j}) }.
\end{equation} In fact, let $N$ denote
$$N:=\#\Set{(s,t)| s<t,\ v^{-1}(s)>v^{-1}(t),\ \check{v}v^{-1} (s)<\check{v}v^{-1} (t) },$$
and let us compute the following numbers $A$ and $B$
\begin{align*}
A&:=\#\Set{(i,j)| i<j,\ v^{-1}(i)>v^{-1}(j),\ \check{v}v^{-1} (i)>\check{v}v^{-1} (j) }\\
&=\#\Set{(i,j)| i<j,\ v^{-1}(i)>v^{-1}(j)}-\#\Set{(i,j)| i<j,\ v^{-1}(i)>v^{-1}(j),\ \check{v}v^{-1} (i)<\check{v}v^{-1} (j) }\\
&=\ell_d{(v)}-N.\\
\\
B&:=\#\Set{(i,j)| i<j,\ v^{-1}(i)<v^{-1}(j),\ \check{v}v^{-1} (i)>\check{v}v^{-1} (j) }\\
&=\#\Set{(i,j)| v^{-1}(i)<v^{-1}(j),\ \check{v}v^{-1} (i)>\check{v}v^{-1} (j) }-\#\Set{(i,j)| i>j,\ v^{-1}(i)<v^{-1}(j),\ \check{v}v^{-1} (i)>\check{v}v^{-1} (j) }\\
&=\#\Set{l<k | \check{v}(l)>\check{v}(k) }-\#\Set{(i,j)| i>j,\ v^{-1}(i)<v^{-1}(j),\ \check{v}v^{-1} (i)>\check{v}v^{-1} (j) }\\
&=\ell_d(\check{v})-\#\Set{(i,j)| i>j,\ v^{-1}(i)<v^{-1}(j),\ \check{v}v^{-1} (i)>\check{v}v^{-1} (j) }\\
&=\ell_d(v)-N.
\end{align*} 
We therefore have $$\ell_d(\check{v}v^{-1}) =\#\Set{1\leq i<j\leq d | \check{v}v^{-1} (i)>\check{v}v^{-1} (j) }=A+B=2\ell_d(v)-2N.$$ But \begin{align*}
N&=\#\Set{(s,t)| s<t,\ v^{-1}(s)>v^{-1}(t),\ \check{v}v^{-1} (s)<\check{v}v^{-1} (t) }\\
&=\#\Set{(s,t)| s<t,\ \overline{v^{-1}(s)}<\overline{v^{-1}(t)},\ \overline{v \big{(}\overline{v^{-1} (s)}\big{)}}<\overline{v \big{(}\overline{v^{-1} (t)}}\big{)}}\\
&=\#\Set{(s,t)| s<t,\ \overline{v^{-1}(s)}<\overline{v^{-1}(t)},\ {v \big{(}\overline{v^{-1} (s)}\big{)}}>{v \big{(}\overline{v^{-1} (t)}}\big{)}}\\
&=\#\set{1\leq i<j\leq d | v (i)>v (j) \text{ and } v (\overline{i})<v (\overline{j}) }.
\end{align*} Whence the equality (\ref{eq:longueur sigmacheck sigmamoinsun}).
\item[$(iii)$] We have the following facts
$$\begin{matrix}\mbox{If } d=2k\mbox{, there exists } \sigma\in \mathfrak{S}_{d}\mbox{ such that }\check{\sigma}=\sigma \mbox{ and}\\ 
v\sigma (i)=min\set{ v\sigma (j)| \ j\in [i,\overline{i}] }\mbox{ for all }i\in \lbrace 1,...,k \rbrace.
\end{matrix}$$

$$\begin{matrix}
\mbox{If } d=2k+1\mbox{, there exists } \sigma\in \mathfrak{S}_{d}\mbox{ such that }\check{\sigma}=\sigma \mbox{ and}\\  
v\sigma (i)=min\set{ v\sigma (j)| \ j\in [i,\overline{i}]\setminus \set{k+1}} \mbox{ for all }i\in \lbrace 1,...,k \rbrace.
\end{matrix}$$

They imply the proposition. Indeed, let us assume that $\sigma\in \mathfrak{S}_d$ satisfies these properties and fix $1\leq i<j\leq d$. If $d=2k$, then if $i\leq k$, we deduce from the minimality property of $v\sigma$ that $v\sigma (i)<v\sigma (j)$ when $j<\overline{\imath}$, and that $v\sigma (\bar{j})<v\sigma (\overline{\imath})$ when $\overline{\imath}<j$. And if $k<i$ then $\bar{j}<i$, $k\leq \overline{\imath}$ and we have again $v\sigma (\bar{j})<v\sigma (\overline{\imath})$ by minimality. In the case where $d=2k+1$ is odd the same arguments give the alternative for $i\neq k+1$ and $j\neq k+1$. But if $i=k+1$ then $\bar{j}<k+1<j$ and $v\sigma (\bar{j})<v\sigma (j)$ by minimality. Since $\check{\sigma}=\sigma$, we have $\sigma(i)=\sigma(k+1)=k+1=\sigma(\overline{\imath})$ hence $v\sigma(i)=v(k+1)=v\sigma(\overline{\imath})$. Comparing $v\sigma(i)$ and $v\sigma (j)$ gives then the alternative thanks to the previous inequality. If $j=k+1$ we apply the same argument with $\overline{\imath}$ in the place of $j$.
 
We now describe a suitable $\sigma$ to prove the previous facts. Let first $a$ be the composition of the transpositions $(i\ \overline{i})$ for $1\leq i\leq k$ such that $v(i)>v(\overline{i})$. Let then $b$ be the unique element of $\mathfrak{S}_d$ such that $\check{b}=b$ and that $vab$ increases on $\set{1,...,k}$. The desired permutation is $\sigma:=ab$. 
\end{enumerate}
\end{proof}
\end{prop}

We are now able to state and prove the desired propositions, similar to Proposition \ref{prop:type A T fixed point}, \ref{prop:type A corriger par z}, \ref{prop:type A a imply hypo 2 }, \ref{prop:type A a b quasi imply hypo 1}. We begin with a parametrization. Let $w_0:=\sigma_{0,r}$ and let $\mathcal{I}_r^{1}$ (respectively  $\mathcal{I}_r^{-1}$) denote the set of involutions in $\mathfrak{S}_r$ (respectively the set of involutions without fixed point). For $u\in W_r$, we also introduce the following subset of $G_r/B_r$:
$$\O(u):=\Set{xB_r| I_{-\e,r}{}^{\delta}xI_{-\e,r}x\in {B_rw_0uB_r}}.$$

\begin{prop}\label{prop:types BD T fixed point}
We have the following parametrization of the $Z$-orbits in $G/B$
\[\begin{array}{ccccl}
& & \mathcal{I}_r^{-\e}\times W^P & \to & Z \backslash (G/B)\\ 
& & (u,v) & \mapsto & {Z}\cdot \begin{psmallmatrix} x & &\\ 
& 1 & \\ & & {}^{\d}x^{-1} \end{psmallmatrix}v^{-1}B\mathrm{\ with\ any\ }xB_r\in \O(u).\end{array}\]
As a consequence, the number of $Z$-orbits is finite. Besides, if $\e=1$ (types B, D), all $Z$-orbits contain a $T$-fixed point. If $\e=-1$ (type C), this is the case exactly for the $Z$-orbits parametrized by $(u,v)$ with $u\in \mathcal{I}_r^{-1}$. 
\begin{proof}
In the literature, there exist several parametrizations of the $G_{-\e,r}$-orbits of the flag variety $G_r/B_r$; we will use the one presented in \cite[Proposition 4]{FressePenkov2017}\footnote{Their article is written for the base field $\mathbb{C}$, but the result used here is valid for any base field of any characteristic, see the proof in \cite[Section 3.5]{FressePenkov2017}.} for $\e=1$ and the one presented in \cite[Example 10.3]{RichardsonSpringer1990} for $\e=-1$, $char(\k)\neq 2$. Considering that $x\mapsto I_{-\e,r}{}^{\delta}x^{-1}I_{-\e,r}$ is an involution of $G_r$ whose fixed-point set consists in the symplectic group (if $\e=1$) or the (usual) special orthogonal group (if $\e=-1$, with $char(\k)\neq 2$), these two works ensure the existence of a bijection
$$\begin{array}{ccccl}
& & \mathcal{I}_r^{-\e}& \to &  G_{-\e,r}\backslash(G_r/B_r)\\
& & u & \mapsto & \O(u).
\end{array}.$$ But if we consider the block shapes of $L$ and $L_Z$, we have a surjection
\begin{equation}\label{eq:surjection equivariant LBL}
L/B_L\twoheadrightarrow G_r/B_r
\end{equation}
which is $L_Z\twoheadrightarrow G_{-\e,r}$ and $T \twoheadrightarrow T_r$ equivariant and which induces the bijection \begin{equation}\label{eq:bijection orbites Gre LZ}
\begin{array}{ccccl}
& & G_{-\e,r}\backslash (G_r/B_r)& \to &  L_Z\backslash \left(L/B_L\right)\\
& & G_{-\e,r}\cdot xB_r & \mapsto & L_{Z}\cdot \begin{psmallmatrix} x & &\\ 
& 1 & \\ 
& & {}^{\delta}x^{-1} \end{psmallmatrix}B_{L}.
\end{array}
\end{equation} between the orbit sets. We deduce the bijection 
$$\begin{array}{ccccl}
& & \mathcal{I}_r^{-\e}& \to &  L_Z\backslash \left(L/B_L\right)\\
& & u & \mapsto & L_Z(u):=L_{Z}\cdot \begin{psmallmatrix} x & &\\
 & 1 & \\ 
 & & {}^{\delta}x^{-1} \end{psmallmatrix}B_{L} \mbox{ with any } xB_r\in \O(u).
\end{array}$$
As in type A, we can then apply \cite[Theorem 7.2 (a)]{ChaputFresse}, and we finally have the bijection
$$\begin{array}{ccccl}
& & \mathcal{I}_r^{-\e}\times W^P & \to & Z \backslash (G/B)\\ 
& & (u,v) & \mapsto & Z(u,v):={Z}\cdot \begin{psmallmatrix} x & &\\ 
& 1 & \\ & & {}^{\d}x^{-1} \end{psmallmatrix}v^{-1}B\mbox{ with any } xB_r \in \O(u).\end{array}$$ 
which is the desired parametrization.  

Let us now conclude about the claim on $T$-fixed points. As usual a superscript will denote the set of fixed points for the involved action. Let $(u,v)$ be a couple of parameters. Since $w_0$ is an involution without fixed point (recall $r$ is even for $\e=1$) and $I_{-\e,r}{}^{\delta}\s I_{-\e,r}=w_0\s^{-1} w_0$ for any $\s \in W_r$, it is not difficult to show that
$$\O(u)^{T_r}\neq \emptyset \Leftrightarrow u\mbox{ is conjugate to } w_0 \mbox{ in }W_r\Leftrightarrow u\in \mathcal{I}_r^{-1}.$$
Then it suffices to show that
\begin{equation}\label{eq:equivalence T fixed points}
Z(u,v)^T\neq \emptyset \Leftrightarrow \O(u)^{T_r}\neq \emptyset.
\end{equation}

But, thanks to the previous reference, we have a ($L$ and thus) $T$-equivariant isomorphism (\cite[Proposition 6.4]{ChaputFresse})
$$L/B_L\simeq \left(P\cdot v^{-1}B\right)^{\tau}$$ where $\tau$ is a suitable cocharacter $\mathbb{G}_{m}\to T$. Moreover, its restriction induces an isomorphism (\cite[Theorem 7.2(b)]{ChaputFresse}) $$L_Z(u)\simeq Z(u,v)^{\tau}.$$ Hence we get an isomorphism $$L_Z(u)^{T}\simeq \left(Z(u,v)^{\tau}\right)^{T}=Z(u,v)^T.$$ Moreover, the equivariant surjection (\ref{eq:surjection equivariant LBL}) which induces the bijection (\ref{eq:bijection orbites Gre LZ}) gives an isomorphism $$L_Z(u)^{T}\simeq \O(u)^{T_r}.$$ Composing the above isomorphisms, we finally get
$$Z(u,v)^T\simeq \O(u)^{T_r}$$ which leads to the desired equivalence (\ref{eq:equivalence T fixed points}).
\end{proof}
\end{prop}

\begin{prop} \label{prop:types BCD corriger par z }
Let $v\in W$. Then, there exists $z_0\in Z\cap N_G(T)$ such that $w:=\zeta_{\e,n} (z_0)^{-1}v$ satisfies:
\begin{enumerate}
\item[$(a)$] $s_v=s_w$ and $w^{-1}$ induces $s_wu_0$ on $\set{r+1,...,n-r}$,
\item[$(b)$] $w^{-1} (i)<w^{-1} (j)$ or $w^{-1} (r-i+1)>w^{-1}(r-j+1)$ for any $1\leq i<j\leq r$.
\end{enumerate}
\begin{proof}
Let $v\in W$. The inverse $v^{-1}$ induces on $\set{r+1,...,n-r}$ an element $u\in \mathfrak{S}_{n-2r}$ such that $s_vu$ is in $W_{\e,n-2r}$ (Proposition \ref{prop:fait sv}) so that $s_vu\s_1=u_0$ for a suitable $\s_1\in W_{\e,n-2r}$. Then $v^{-1}\begin{psmallmatrix} 1 & & \\ & \s_1 & \\ & & 1\end{psmallmatrix}$ induces $s_vu_0$ on $\set{r+1,...,n-r}$. Besides, there exists a monomial matrix $g_1\in {G_{\e,n-2r}}$ such that $\zeta_{\e,n-2r}(g_1)=\s_1$. Applying Proposition \ref{prop:préliminaire technique types B,C,D} $(iii)$ on $v^{-1}$ gives also $\s_2\in \mathfrak{S}_r$ such that $v^{-1}\begin{psmallmatrix}
\s_2 & &\\
& I_{n-2r} &\\
& & \check{\s_2}
\end{psmallmatrix}$ satisfies $(b)$ with $\check{\s_2}=\s_2$ relatively to $\set{1,...r}$. We can then find a monomial matrix $g_2\in {G_{-\e,r}}$ such that $\zeta_{-\e,r}(g_2)=\s_2$. Indeed, in type C ($\e=-1$, odd $r$, $char(\k)\neq 2$), we can define $g_2$ as the matrix product $\sigma_2\begin{psmallmatrix}
1 & &\\
 & det(\s_2) & \\
 & & 1
\end{psmallmatrix}$. In types B, D ($\e=1$ and even $r$), since $\check{\s_2}=\s$ in $\mathfrak{S}_r$, there exists $\lambda \in T_{r/2}$, with coefficients $\pm 1$, such that $I_{-1,r}\s_2I_{-1,r}=\begin{psmallmatrix}
\lambda & \\
& {}^{\delta}\lambda
\end{psmallmatrix}\s_2$. We can then define $g_2$ as the matrix product $\begin{psmallmatrix}
\lambda & \\
& 1
\end{psmallmatrix}\s_2$. Now we can introduce
$$z_0:=\begin{pmatrix}
g_2 & &\\
& g_1 &\\
& & I_{-\e,r} g_2 I_{-\e,r}
\end{pmatrix},$$
which is in $Z\cap N_G(T)$ with
$$\zeta_{\e,n}(z_0)=\begin{psmallmatrix}
\zeta_{-\e,r}(g_2) & &\\
& \zeta_{\e,n-2r}(g_1) &\\
& & \zeta_{-\e,r}(I_{-\e,r}g_2I_{-\e,r})
\end{psmallmatrix}
=\begin{psmallmatrix}
\s_2 & &\\
& \s_1 &\\
& & \check{\s_2}
\end{psmallmatrix}.$$
We see that $v^{-1}\zeta_{\e,n}(z_0)$ and $v^{-1}\begin{psmallmatrix} 1 & & \\ & \s_1 & \\ & & 1\end{psmallmatrix}$ induce the same element $s_vu_0$ on $\set{r+1,...,n-r}$. Since $\s_1\in W_{\e,n-2r}$, the quantities $$\#\Set{i\in \set{r+1,...,n-r}| v^{-1}(i)>\lfloor n/2 \rfloor}$$ and $$\#\Set{i\in \set{r+1,...,n-r}| v^{-1}\begin{psmallmatrix} 1 & & \\ & \s_1 & \\ & & 1\end{psmallmatrix}(i)>\lfloor n/2 \rfloor}$$ have the same parity for $\e=1$, even $n$. Hence $s_v=s_{\zeta_{\e,n}(z_0)^{-1}v}$. On the other hand, $v^{-1}\zeta_{\e,n}(z_0)$ and $v^{-1}\begin{psmallmatrix}
\s_2 & &\\
& I_{n-2r} &\\
& & \check{\s_2}
\end{psmallmatrix}$ induce the same element on $\set{1,...,r}$. We conclude that $w:=\zeta_{\e,n}(z_0)^{-1}v$ satisfies the desired conditions.
\end{proof}
\end{prop}

\begin{prop}\label{prop:type BCD a imply hypo 2}
Suppose $w\in W$ satisfies the property $(a)$ of Proposition \ref{prop:types BCD corriger par z }. Then $z^{-1}\varpi (z)\in \overline{\conj{w}{B}B}$ for all $z\in Z$.
\begin{proof}
Let $w^{-1}$ be such an element and fix $$w_0:=\begin{psmallmatrix}
1 & & \\
 & u_0s_w & \\
 & & 1
\end{psmallmatrix}.$$ We have
\begin{equation} \label{eq:1 condition i theo type BCD}
\begin{pmatrix}
1 &  & \\
 & ^{u_0} B_{\e,n-2r} & \\
 &  & 1
\end{pmatrix} \subset \ {}^wB.
\end{equation}
Indeed, ${s_w}$ stabilizes $B_{\e,n-2r}$ because, in the case $s_w\neq id$, even $n-2r$, we have for all $i,j\in \set{1,...,n-2r}$ with $1\leq i<j<\overline{i}=n-2r-i+1$ the inequality $i<(n-2r)/2$, so that $s_w(i)=i<s_w(j)$. The inclusion (\ref{eq:1 condition i theo type BCD}) is thus equivalent to:
$${w^{-1}w_0}\begin{pmatrix}
1 &  & \\
 & B_{\e,n-2r} & 0\\
 &  & 1
\end{pmatrix}{w_0^{-1}w} \subset B,$$
that is  $$\forall r+1 \leq i<j\leq n-r ,\ w^{-1}w_0(i)<w^{-1}w_0(j).$$
But, $w^{-1}$ induces $s_wu_0$ on $\set{r+1,...,n-r}$, which means that $$w^{-1}\begin{pmatrix}
1 & & \\
 & (s_wu_0)^{-1} & \\
 & & 1
\end{pmatrix}=w^{-1}w_0$$ increases on $\set{r+1,...,n-r}$. We get (\ref{eq:1 condition i theo type BCD}). But $u_0$ is the longest element in $W_{\e,n-2r}$ so that $${}^{u_0}B_{\e,n-2r}B_{\e,n-2r}\simeq B_{\e,n-2r}u_0B_{\e,n-2r}$$ is dense in $G_{\e,n-2r}$ so that (\ref{eq:1 condition i theo type BCD}) implies 

\begin{equation} \label{eq:2 condition i theo type BCD}
\begin{pmatrix}
1 & & \\
 & G_{\e,n-2r} & \\
 & & 1
\end{pmatrix}\subset \overline{{}^wBB}.
\end{equation}

Now let $z\in Z$. The element $z^{-1}\varpi (z)$ has the shape 
$$\begin{pmatrix}
1 & A_1 & A_2\\
0 & C & A_3\\
0 & 0 &1 
\end{pmatrix}=\begin{pmatrix}
1 & & \\
 & C & \\
 & &1 
\end{pmatrix} \begin{pmatrix}
1 & A_1 & A_2\\
0 & 1 & C^{-1}A_3\\
0 & 0 &1 
\end{pmatrix}$$ with $C\in G_{\e,n-2r}$ and suitable matrices $A_i$. We conclude $z^{-1}\varpi (z)\in \overline{^wBB} B=\overline{^wBB}$ by (\ref{eq:2 condition i theo type BCD}). 
\end{proof}
\end{prop}

\begin{prop}\label{prop:types BCD a b quasiimply hypo 1}
Suppose $w\in W$ satisfies properties $(a)$ and $(b)$ of Proposition \ref{prop:types BCD corriger par z }. Let $w=\tau\nu$ be the decomposition of $w$ in $W_P(W^P)^{-1}$. Then $$\ell (\tau) -\ell (\tau ^{-1}\theta (\tau))/2=\dim Z/Z\cap B-\dim H/B_H.$$
\begin{proof}
Let $w$ be such an element. There exists $\s \in \mathfrak{S}_{r}$ such that $w^{-1}\s$ is increasing on $\set{1,...,r}$. We fix then
$$\tau:= \begin{psmallmatrix}
\sigma & 0 & 0\\
0 & u_0 & 0\\
0 & 0 & \check{\sigma}
\end{psmallmatrix}$$ which is in $W_P$. With the decomposition $$w^{-1}\tau=w^{-1}\begin{psmallmatrix} 1 & & \\ & (s_wu_0)^{-1} & \\ & & 1 \end{psmallmatrix} \begin{psmallmatrix} \s & & \\ & s_w & \\ & & \check{\s}\end{psmallmatrix}$$ we see that $w^{-1}\tau$ is increasing on $\set{1,...,r}$, on $\set{r+1,...,n-r}$ and that 
$${w^{-1}\tau(\lfloor n/2 \rfloor)}<{w^{-1}\tau(\lfloor n/2 \rfloor+1+\frac{\e+(-1)^n}{2})}.$$ Thus, $w^{-1}\tau$ is in $W^P$ and $w^{-1}\tau=v^{-1}$. Applying $\theta$, we find
$$ \tau ^{-1}\theta (\tau)=\begin{psmallmatrix}
\sigma ^{-1 }& & \\
 & u_0^{-1} & \\
& & \check{\sigma}^{-1} 
\end{psmallmatrix} \begin{psmallmatrix}
\check{\sigma} & & \\
 & u_0 & \\
 & & {\sigma}
\end{psmallmatrix}
= \begin{psmallmatrix}
\sigma^{-1} \check{\sigma}& & \\
 & 1 & \\
 & & \check{\sigma}^{-1} {\sigma}
\end{psmallmatrix}$$ so that $\ell_n( \tau ^{-1}\theta (\tau))=2\ell_{r}(\sigma ^{-1}\check{\sigma})$. Besides $\ell_n( \tau)=\ell_{n-2r}(u_0)+2\ell_r(\sigma)$. Using Proposition \ref{prop:préliminaire technique types B,C,D} $(i)$ we see that $\ell_n( \tau ^{-1}\theta (\tau))=2\ell( \tau ^{-1}\theta (\tau))$ and $\ell_n(\tau)-2\ell(\tau)=\ell_{n-2r}(u_0)-2\ell(u_0)$. Since $\begin{psmallmatrix}
\s & &\\
&  I_{n-2r}&\\
& & \check{\s}
\end{psmallmatrix}=\left(w^{-1}\begin{psmallmatrix}
\s & &\\
&  I_{n-2r}&\\
& & \check{\s}
\end{psmallmatrix}\right)^{-1}w^{-1}$ where $\left(w^{-1}\begin{psmallmatrix}
\s & &\\
&  I_{n-2r}&\\
& & \check{\s}
\end{psmallmatrix}\right)^{-1}$ increases on $w^{-1}(\set{1,...,r})$, $\sigma$ satisfies property $(b)$ of Proposition \ref{prop:types BCD corriger par z }. By Proposition \ref{prop:préliminaire technique types B,C,D} $(ii)$ we have $\ell_r(\check{\sigma}^{-1}\sigma)=2\ell_r(\sigma^{-1})$ that is $\ell_r(\sigma^{-1}\check{\sigma})=2\ell_r(\sigma)$. Therefore we have $\ell(\tau ^{-1}\theta (\tau))/2=\ell_n(\tau ^{-1}\theta (\tau))/4=\ell_r(\sigma)$ and $\ell(\tau)=1/2\big{(}\ell_n(\tau)-\ell_{n-2r}(u_0)+2\ell(u_0)\big{)}=\ell_r(\sigma)+\ell(u_0)$. Hence $\ell(\tau)-\ell(\tau ^{-1}\theta (\tau))/2=\ell(u_0)$. An easy computation gives $\dim Z/Z\cap B-\dim H/B_H=\dim G_{\e,n-2r}/B_{\e,n-2r}=\ell(u_0)$ and we have the desired formula.
\end{proof}
\end{prop}

\subsubsection{Conclusion}\label{sec:conclusion birationality}
We can now prove the claims of Theorem \ref{theo:1 birational rational} with the exception of rationality. We recall that $G$ is the group $Gl_{n\k}$, $SO_{n\k}$, or $Sp_{n\k}$ under certain restrictions (see Section \ref{sec:application through matrix models} for the full context). So let $Y$ be a $Z$-orbit closure in $G/B$.  In type C, we assume that the dense $Z$-orbit possesses a $T$-fixed point. By hypothesis (type C), Proposition \ref{prop:type A T fixed point} (type A) or Proposition \ref{prop:types BD T fixed point} (types B, D), there exists $w\in W$ such that $Y=\overline{Z\cdot wB}$. By combining Propositions \ref{prop:type A corriger par z} and \ref{prop:type A a imply hypo 2 } (type A) or Propositions \ref{prop:types BCD corriger par z } and \ref{prop:type BCD a imply hypo 2} (types B, C, D), we can assume without loss of generality that $w$ satisfies the hypothesis $2$ of Lemma \ref{lem:équivalence birationnelle cas général}. The hypothesis $4$ is satisfied by Proposition \ref{prop:smoothness ZcapwB} whereas the hypothesis $3$ is satisfied by the point $(i)$ of Lemma \ref{lem:formule dimension sous-groupe symétrique} with Proposition \ref{prop:connexity}. For $char(\k)\neq 2$, the point $(ii)$ of the same Lemma together with Proposition \ref{prop:type A a b quasi imply hypo 1} (type A) or \ref{prop:types BCD a b quasiimply hypo 1} (types B, C, D) ensure that the hypothesis $1$ is satisfied. But, thanks to the smoothness given by Proposition \ref{prop:smoothness ZcapwB}, the same dimension formula holds for $char(\k)=2$ and types A, B, D. We can thus apply Lemma \ref{lem:équivalence birationnelle cas général} and get the desired result. 

\section{Normality, rationality, Cohen-Macaulayness, general case}\label{sec:normality rationality}

We prove a stronger version of Theorem \ref{theo:2 rationality normality general argument}. We first clarify the notation used in this statement and the general context. 

\subsection{Context}\label{sec:context general types free}

Let $G$ be a connected semi-simple algebraic group over $\mathds{k}$ and $(T,B)$ a Killing pair. Let $H\subset G$ be a connected reductive subgroup such that $(T_H,B_H):=(T\cap H, B\cap H)$ is a Killing pair. For any character $\lambda$ of $T$, let $\mathds{k}_{\lambda}$ be the one dimensional representation of $B$ with weight $\lambda$ and $\mathcal{L}_G(\lambda)$ be the $G$-equivariant line bundle on $G/B$ corresponding to $G\times^B\mathds{k}_{-\lambda}\to G/B$. Let $$V_G (\lambda):=H^0(G/B,\mathcal{L}_G(\lambda))$$ denotes the dual Weyl $G$-module with lowest weight $-\lambda$ (if $\lambda$ is dominant).

We fix $w$ in the Weyl group $W$ of $G$ and consider the Schubert variety $\overline{B\cdot wB}\subset G/B$. We also consider the natural morphisms
$$q\colon H\times ^{B_H} \overline{B\cdot wB}\to G/B, [h,gB]\mapsto hgB$$ and $$k\colon H\times ^{B_H} \overline{B\cdot wB}\to H/B_H, [h,gB]\mapsto hB_H.$$ 
Let $Z_{\mathfrak{w}}$ be the Bott-Samelson variety associated with the choice of a reduced word $\mathfrak{w}$ that decomposes $w$ into simple reflections. The $B_H$-equivariant Bott-Samelson resolution $Z_{\mathfrak{w}}\to \overline{B\cdot wB}$ induces a birational morphism $H\times^{B_H}Z_{\mathfrak{w}}\to H\times^{B_H}\overline{B\cdot wB}$. By composing with $q$ and $k$, we obtain morphisms
$$\overset{\sim}{q}\colon H\times ^{B_H} Z_{\mathfrak{w}}\to G/B,$$
$$\overset{\sim}{k}\colon H\times ^{B_H}Z_{\mathfrak{w}} \to H/B.$$
Let $\pi$ and $\overset{\sim}{\pi}$ be the restrictions of $q$ and $\overset{\sim}{q}$ onto their respective images 
$$\pi\colon H\times ^{B_H} \overline{B\cdot wB}\to\overline{HB\cdot wB},$$
 $$\overset{\sim}{\pi}\colon H\times ^{B_H} Z_{\mathfrak{w}}\to\overline{HB\cdot wB}.$$

\begin{theo}\label{theo:theorem 2 plus}
We make the following assumptions:
\begin{enumerate}
\item[(i)] the morphism $\pi\colon H\times ^{B_H} \overline{B\cdot wB}\to \overline{HB\cdot wB}$ is birational,
\item[(ii)] the character $2\rho_H-\rho_{G\vert T_H}$ is dominant,
\item[(iii)] $char(\k)=0$ or 
\item[(iii)'] $char(\k)=p>0$ and the restriction $V_G\left((p-1)\rho_G\right)\to V_H\left((p-1)\rho_{\vert T_H}\right)$ is surjective.
\end{enumerate}
Then $\overline{HB\cdot wB}$ is normal and $\overset{\sim}{\pi}$ is rational. 
Moreover if
\begin{enumerate}
\item[(iv)]  $char(\k)=0$ or
\item[(iv)'] there exists a line bundle $\mathcal{M}$ on $G/B$ such that $k^*\mathcal{L}_H(\rho_{G\vert T_H}-2\rho_H)\simeq q^*\mathcal{M}$,
\end{enumerate} 
then $\overline{HB\cdot wB}$ is Cohen-Macaulay with dualizing sheaf $\overset{\sim}{\pi}_*\w_{H\times ^{B_H} {Z_{\mathfrak{w}}}}$ and we have the vanishings $R^i\overset{\sim}{\pi}_*\w_{H\times ^{B_H} {Z_{\mathfrak{w}}}}=0$ for all $i>0$.
\end{theo}

\subsection{Normality and rationality}\label{sec:normality and rationality}

We start with the proof of the first part of Theorem \ref{theo:theorem 2 plus} concerning normality and rationality.

\subsubsection{An inductive result by Perrin and Smirnov}

We present here a general setup which contains as a particular case the Bott-Samelson resolution. Let $Y$ be a scheme and $n$ an integer. Let us consider, for $i\in [0,n]$, the following schemes and morphisms:
\begin{enumerate}
\item[-] a scheme $T_i$ and a morphism $Y\overset{p_i}{\longrightarrow} T_i$,
\item[-] schemes $\overset{\sim}{X}_i$ and $X_i$ over $Y$,
\item[-] morphism $\overset{\sim}{X}_i\to X_i$ over $Y$,
\end{enumerate}
such that, for all $i\in [0,n-1]$,
\begin{equation}\label{eq:définition tildeXi}
\xymatrix{\overset{\sim}{X}_{i+1} \ar[r] \ar[dd] & \overset{\sim}{X}_i \ar[d]\\
& Y\ar[d]^{p_{i+1}}\\
Y\ar[r]_{p_{i+1}} & T_{i+1}}
\end{equation}
is cartesian and for all $i\in [0,n]$, $X_i$ is the scheme-theoretic image of $\overset{\sim}{X}_i\to Y$:
\begin{equation}\xymatrix{\overset{\sim}{X}_{i} \ar[rr] \ar[dr] & & Y.\\
& X_i\ar@{^{(}->}[ur]}
\end{equation}
Besides, we will say that a morphism $f:Z\to T$ satisfies (\ref{eq:properties on pi}) if the following hold:
\begin{equation}\label{eq:properties on pi}
\begin{split}
& f\colon  Z\to T \mbox{ is faithfully flat and proper, }\\
&\mbox{its geometric fibers } Z_{\overline{t}} \mbox{ are  connected, normal and reduced,}\\
&\mbox{with }\dim Z_t\leq 1,\ H^1(Z_t,\O_{Z_t})=0.  
\end{split}\tag{*}
\end{equation}

Let also \textbf{P} be a property of morphism of schemes which is preserved under any composition, any base change and which is satisfied by closed immersions. 

The article \cite[]{PerrinSmirnov} by Perrin and Smirnov then contains the following result, inspired by \cite[]{BrionKumar}. The proof consists of increasing and decreasing induction based on appropriate cartesian diagrams and the stability under base change, descent, and cancellation for certain properties of morphisms of schemes. The authors state this result for varieties, but we present a statement in terms of schemes to emphasize the generality of their argument.

\begin{prop}\label{prop:PerrinSmirnov induction}
We assume that for all $i$, $p_i$ satisfy (\ref{eq:properties on pi}), $T_i$ is locally noetherian and there exists an ample line bundle $\mathcal{M}_i$ on $T_{i}$ such that the restriction 
 $$H^0\big{(}X_{i+1}, (p_{i+1}^*\mathcal{M}_{i+1})_{\vert X_{i+1}}\big{)} \to H^0\big{(}X_{i}, (p_{i+1}^*\mathcal{M}_{i+1})_{\vert X_{i}}\big{)}$$ 
 is surjective. We also make the following assumptions: 
 \begin{enumerate}
     \item[-] the morphism $\overset{\sim}{X}_0\to X_0$ possesses \textbf{P} and is rational,
     \item[-] the morphism $\overset{\sim}{X}_0\to X_0$ is birational,
     \item[-] the scheme $X_0$ is normal.
 \end{enumerate}
 Then for all $i$, the morphism $\overset{\sim}{X}_i\to X_i$ possesses $\textbf{P}$, is rational, birational and the scheme $X_i \textit{ is normal}$.
\end{prop}
We will apply this result to suitable families. Recall the data $G$, $H$, $B$, $T$, $w$, ...  that we have fixed and let $n:=\ell(w)$. Let $P_{\a_k}$ be the minimal parabolic subgroup relative to the simple root $\a_k$. Assume that we have $\mathfrak{w}=(s_{\a_{j_1}},...,s_{\a_{j_n}})$. For all $i\in [0,n]$ let $\mathfrak{w}_i$ denote the subword $(s_{\a_{j_1}},...,s_{\a_{j_i}})$ and $w_i:=s_{\a_{j_1}}...s_{\a_{j_i}}$ the corresponding Weyl group element. Let $Z_{\mathfrak{w}_i}$ be the Bott-Samelson variety associated to $\mathfrak{w}_i$. As we did previously for $\mathfrak{w}=\mathfrak{w}_n$ with $q$ and $\overset{\sim}{q}$, we define $q_i$ as the composition \begin{equation}
q_i\colon H\times^{B_H} \mathcal{Z}_{\mathfrak{w}_i} \to G/B.
\end{equation} We then denote: \begin{align*}
Y&:=G/B\\
T_{i}&:=G/P_{\a_{j_i}}\\
p_i&:=G/B\to G/P_{\a_{j_i}}\\
\overset{\sim}{X}_i&:=H\times ^{B_H} \mathcal{Z}_{\mathfrak{w}_i},\\
\overset{\sim}{X}_i\to Y&:= q_i,\\
X_i&:=Im\ q_i=\overline{HB\cdot w_iB},
\end{align*} and we choose an ample line bundle $\mathcal{M}_i$ on the projective variety $T_i$. By construction of the Bott-Samelson variety, we have:
\begin{align*}
\overset{\sim}{X}_i\times _{T_{i+1}}Y&=\left(H\times^{B_H} \mathcal{Z}_{\mathfrak{w}_i}\right)\times_{G/P_{\a_{j_{i+1}}}} G/B\\
&\simeq H\times^{B_H} \left( \mathcal{Z}_{\mathfrak{w}_i}\times_{G/P_{\a_{j_{i+1}}}}G/B\right)\\
&\simeq H\times^{B_H} \mathcal{Z}_{\mathfrak{w}_{i+1}}\\
&=\overset{\sim}{X}_{i+1}.
\end{align*} Besides, we note that the $P_{\a_{j_i}}/B\simeq \mathbb{P}^1$-fibration $p_i$ satisfies (\ref{eq:properties on pi}) and that $$\overset{\sim}{X_0}\simeq X_0\simeq H/B_H$$ are normal varieties.  On the other hand, the hypothesis $(i)$ on the birationality of $\pi$ ensures by composition that $\overset{\sim}{\pi}$ is birational, that is: $$\overset{\sim}{X}_n\to X_n \text{ is birational}.$$ Because the $X_i$ are embedded into $Y$ in a way compatible with their mutual inclusions and the pullbacks $p_{i}^*\mathcal{M}_i$ are semi-ample on $Y=G/B$, it is now sufficient to prove the surjectivity
\begin{equation}\label{eq:surjectivity of restriction}
H^0\left(G/B, \L\right)\twoheadrightarrow H^0\left(\overline{HB\cdot w_iB}, \L_{\vert \overline{HB\cdot w_iB}}\right)
\end{equation}
for all $i$ and all semi-ample line bundles $\mathcal{L}$ on $G/B$. Once we show (\ref{eq:surjectivity of restriction}), we get the surjectivity
 $$H^0\big{(}X_{i+1}, (p_{i+1}^*\mathcal{M}_{i+1})_{\vert X_{i+1}}\big{)} \twoheadrightarrow H^0\big{(}X_{i}, (p_{i+1}^*\mathcal{M}_{i+1})_{\vert X_{i}}\big{)}$$ 
for all $i$, which allows us to apply the previous proposition so that we obtain the first part of the theorem.

We fix now $i\in \set{0,...,n}$, a semi-ample line bundle $\mathcal{L}$ on $G/B$ and we prove (\ref{eq:surjectivity of restriction}). We will need to distinguish between the zero and positive characteristic cases.

\subsubsection{Positive characteristic case}

We first assume that $p:=char(\mathds{k})>0$. The hypotheses $(ii)$ and $(iii)'$ of the Theorem \ref{theo:theorem 2 plus} enable us to apply the main result of He and Thomsen in \cite[Theorem 20]{HeThomsen}. It gives a Frobenius splitting of $G/B$, relative to the line bundle $\mathcal{L}_G\left((p-1)\rho_G\right)$ and which compatibly splits our variety $\overline{HB\cdot w_iB}$. Because $\mathcal{L}_G\left((p-1)\rho_G\right)$ is ample and $\mathcal{L}$ semi-ample, \cite[Theorem 1.4.8]{BrionKumar} gives: \begin{equation}\label{eq:vanishing and surjectivity positive characteristic}
\begin{split}
H^1\left(G/B, \L\right)&=0,\\
H^0\left(G/B, \L\right)&\twoheadrightarrow H^0\left(\overline{HB\cdot w_iB}, \L_{\vert \overline{HB\cdot w_iB}}\right),
\end{split}
\end{equation}and in particular we have the desired surjectivity (\ref{eq:surjectivity of restriction}).

\subsubsection{Characteristic zero}

We now show how the previous surjectivity result can be extended from positive to zero characteristic. So let us assume in the sequel that $char(\mathds{k})=0$. We start by realizing our data $G$, $B$, $H$, $X_i$, $w_i$, $\mathcal{L}$, ..., as schemes, morphisms or sheaves over a suitable base ring. This can be done using the following lemma, which also preserves the dominance property of the character $2\rho_H-\rho_{G\vert T_H}$ on all geometric fibers. 

\begin{lem}\label{lem:realisation total}
There exists a $\mathbb{Z}$-subalgebra $A$ of $\mathds{k}$ of finite type,  group-schemes $\mathcal{G}$, $\mathcal{H}$, $\mathcal{T}$, $\mathcal{B}$, $\mathcal{W}$ group schemes over $A$ and a section $\smallw_i \in \mathcal{W}(A)$ such that
\begin{enumerate}
\item[-] $\mathcal{G}$ is the semi-simple Chevalley group scheme over $A$ with $\mathcal{G}_{\mathds{k}}=G$, and with $(\T,\B)$ as a Killing pair satisfying $\mathcal{T}_{\mathds{k}}=T$, $\mathcal{B}_{\mathds{k}}=B$.
\item[-] $\mathcal{H}$ is a closed subgroup of $\G$ and is the reductive Chevalley group scheme over $A$ with $\mathcal{H}_{\mathds{k}}=H$ and with $(\T_{\H},\B_{\H}):=(\T\times_{\G}\H,\B\times_{\G}\H)$ as a Killing pair satisfying ${\T_{\H}}_{\k}=T_H$, ${\B_{\H}}_{\k}=B_H$.
\item[-] The fppf quotient sheaf $\mathcal{G}/ \mathcal{B}$ is representable by an $A$-scheme of finite presentation, smooth and projective, and such that its base change $(\mathcal{G}/ \mathcal{B})_{K}$ over any algebraically closed field $K$ is the flag variety $ \mathcal{G}_{K}/ \mathcal{B}_{K}$.
\item[-] $\mathcal{W}$ is the Weyl group of $\mathcal{G}$ related to $\mathcal{T}$ and the base change ${\smallw_i}_{\mathds{k}}$ recovers the element $w_i\in W$.
\end{enumerate} 
Moreover, there also exists a line bundle $\mathbb{L}$ over $\G / \B$ and an $A$-scheme $\X _i$ flat and projective, equipped with a closed immersion $\X _i\to \G / \B$ over $A$, such that over $G/B$
 \begin{equation}\label{eq:realization linebundle base change over K}
{\mathbb{L}}_{\mathds{k}}\simeq \mathcal{L},
\end{equation}
 \begin{equation}\label{eq:realization our varieties base change over K}
(\mathcal{X}_i)_{\mathds{k}}\simeq
\overline{HB\cdot w_iB}
\end{equation}
and for all $s\in \mathrm{\mathrm{Spec}}\ A$, 
\begin{equation}\label{eq:lem dominance caractere fibres geo}
2\rho_{{\mathcal{H} }_{\overline{s}}}-\rho_{\G_{\overline{s}}\vert {\mathcal{T}_{\mathcal{H} }}_{\overline{s}}} \mathrm{\ is\ dominant}
\end{equation} and we have over $(\mathcal{G}/{\mathcal{B}})_{\overline{s}}$ {\begin{equation}\label{eq:realization our varieties base change}
{\mathcal{X}_i}_{\overline{s}}\simeq \overline{ {\mathcal{H} }_{\overline{s}}{\mathcal{B} }_{\overline{s}}\cdot {{\smallw_i} }_{\overline{s}}{\mathcal{B} }_{\overline{s}}}.
\end{equation}}

\end{lem}
A proof is in the Appendix A. Taking into account classical results on general schemes (\cite[§8]{EGAIV8a15}) and group schemes (\cite[Exposés XIX to XXVI]{SGA3XIXaXXVI}) let us emphasize that our work focuses on the realization of the variety $X_i=\overline{HB\cdot w_iB}$ (see Corollary \ref{cor:réalisation HBwB barre}). It indeed demands caution as it raises a problem of scheme-theoretic image formation under non flat base changes (see Theorem \ref{theo:relèvement image} and Corollary \ref{cor:réalisation H1H2...HnZ barre}).

The base $\mathrm{Spec}A$ is appropriate in the sense that the residue fields $\kappa(x)$ of its points $x$ produce almost all characteristics. 

\begin{lem}\label{prop:car Z algebra finite type}
Let $A$ be an integral $\mathbb{Z}$-algebra of finite type whose characteristic is zero. Then the set $\Set{char(\kappa(x))| x\in \mathrm{\mathrm{Spec}}\ A}$ contains all but finitely many primes. 
\begin{proof}
Since the characteristic is zero, we have $\mathbb{Z}\subset A$ and let $f\colon \mathrm{Spec}\ A\to \mathrm{Spec}\ \mathbb{Z}$ be the finite type dominant morphism related to this inclusion. By Chevalley Theorem $f(\mathrm{Spec}\ A)$ is constructible, and then contains an open dense subset $U$ of its closure $\mathrm{Spec}\ \mathbb{Z}$. Hence, the closed subset $\mathrm{Spec}\ \mathbb{Z}\setminus U$ is finite. On the other hand, we have $char(\kappa(x))=p$ for all prime $p\in \mathbb{Z}$ and $x\in f^{-1}(\set{p\mathbb{Z}})$. We then have $\Set{p| p\text{ prime, } p\mathbb{Z}\in U}\subset\Set{char(\kappa(x))| x \in \mathrm{Spec}\ A}$ and the desired result.
\end{proof}
\end{lem}

We can then choose $s\in \mathrm{Spec}\ A$ such that $p:=char(\kappa(s))$ is large enough to satisfy the assumption of \cite[Lemma 14]{HeThomsen}.  It gives the surjectivity: $$V_{{{\mathcal{G}} }_{\overline{s}}}\big((p-1)\rho_{{{\mathcal{G}} }_{\overline{s}}}\big)\to V_{{{\mathcal{H}} }_{\overline{s}}}\big((p-1){\rho_{{{\mathcal{G}} }_{\overline{s}}}}_{\vert {{{\mathcal{T}}_{\mathcal{H}}}_{\overline{s}}}}\big).$$

We are now able to use the result (\ref{eq:vanishing and surjectivity positive characteristic}) of the previous section on vanishing and surjectivity. Thanks to the isomorphism (\ref{eq:realization our varieties base change}), we have $$H^1 \left(({{\mathcal{G}/{\mathcal{B}}}) }_{\overline{s}}, { \mathbb{L} }_{\overline{s}}\right)=0$$ and the surjectivity
$$H^0\left(({{\mathcal{G}/{\mathcal{B}}}) }_{\overline{s}},{ \mathbb{L} }_{\overline{s}}\right)\twoheadrightarrow H^0\left({\mathcal{X}_{i}}_{\overline{s}},{{ \mathbb{L} }_{\overline{s}}}_{\vert {\mathcal{X}_{i}}_{\overline{s}}}\right).$$ By semi-continuity, the surjectivity can then pass from these geometric special fibers over $s$ to the generic one (\cite[Lemma 1.6.3]{BrionKumar}): we have 
$$H^0\left(({{\mathcal{G}/{\mathcal{B}}}) }_{\overline{\eta}},{ \mathbb{L} }_{\overline{\eta}}\right)\twoheadrightarrow H^0\left({\mathcal{X}_{i}}_{\overline{\eta}},{{ \mathbb{L} }_{\overline{\eta}}}_{\vert {\mathcal{X}_{i}}_{\overline{\eta}}}\right),$$ where $\eta\in \mathrm{Spec}\ A$ denotes the generic point. Tensorizing with the field extension $\overline{\kappa(\eta)}\subset \mathds{k}$ gives then
$$H^0\left({\left({{\mathcal{G}/{\mathcal{B}}}) }_{\overline{\eta}}\right)}_{\mathds{k}},({{ \mathbb{L} }_{\overline{\eta}}})_{\mathds{k}}\right)\twoheadrightarrow H^0\left({({\mathcal{X}_{i}}_{\overline{\eta}})}_{\mathds{k}},{({{ \mathbb{L} }_{\overline{\eta}}})_{\mathds{k}}}_{\vert {({\mathcal{X}_{i}}_{\overline{\eta}})}_{\mathds{k}}}\right),$$ but we have
$${(({{\mathcal{G}/{\mathcal{B}}}) }_{\overline{\eta}})}_{\mathds{k}}\simeq ({\mathcal{G}/{\mathcal{B}}})_{\mathds{k}}\simeq G/B$$ and, with (\ref{eq:realization linebundle base change over K}), $$({{ \mathbb{L} }_{\overline{\eta}}})_{\mathds{k}}\simeq ({ \mathbb{L} })_{\mathds{k}}\simeq  \mathcal{L}$$ and, with (\ref{eq:realization our varieties base change over K}), 
$$({\mathcal{X}_{i}}_{\overline{\eta}})_{\mathds{k}}\simeq {\mathcal{X}_{i}}_{\mathds{k}}\simeq \overline{HB\cdot w_iB}.$$
We therefore recognize the desired surjectivity (\ref{eq:surjectivity of restriction}).

\subsection{Cohen-Macaulayness and rational singularities}

The second part of Theorem \ref{theo:theorem 2 plus} will then follow if we can prove that
\begin{equation}\label{eq:vanishing canonical sheaf}
R^j\overset{\sim}{\pi}_*\w_{\overset{\sim}{X}_n}=0,\ j>0.
\end{equation} In the characteristic zero case, (\ref{eq:vanishing canonical sheaf}) is automatically verified thanks to the Grauert-Riemenschneider theorem \cite[]{GrauertRiemenschneider}. Now let us assume that the characteristic is positive and that $k^*\mathcal{L}_H(\rho_{G\vert T_H}-2\rho_H)\simeq q^*\mathcal{M}$ for a suitable line bundle $\mathcal{M}$ (hypothesis $(iv)'$).

\subsubsection{Another Perrin and Smirnov's argument}

We again use the article \cite[]{PerrinSmirnov} of Perrin and Smirnov. We consider the data involved in Proposition \ref{prop:PerrinSmirnov induction}, to which we add the following. Denoting by $\phi_i$ the morphisms $\overset{\sim}{X}_{i+1}\to \overset{\sim}{X}_i$ we introduce their natural sections $\sigma_i$ by the commutative diagrams
\begin{equation}
\xymatrix{ & & Y \ar[dr]^{p_{i+1}}\\
\overset{\sim}{X}_i \ar[r]^{\s_i} \ar[urr]^{q_i} \ar[drr]_{id} & \overset{\sim}{X}_{i+1}\ar[ur]_{q_{i+1}} \ar[dr]^{\phi_i} & & T_{i+1}.\\
& & \overset{\sim}{X}_i \ar[ur]_{p_{i+1}q_i}}
\end{equation}

We assume that the morphisms $\overset{\sim}{X}_i\to X_i$ are proper birational and that such $p_i$ (and thus the $\phi_i$) is a $\mathbb{P}^1$-fibration over a field. Hence the scheme-theoretic images $Im\ \s_i$ are closed subschemes of codimension one of $\overset{\sim}{X}_{i+1}$ and we can inductively define the divisors 
$$\partial \overset{\sim}{X}_{i+1}=\sigma_i( \overset{\sim}{X}_i)\cup \phi_{i}^{-1}(\partial \overset{\sim}{X}_i).$$ Finally, we put 
$$\Phi_i:=\phi_0\phi_1..\phi_{i-1}.$$

By an argument similar to that of Brion and Kumar in \cite[Proposition 2.2.2]{BrionKumar}\footnote{That is, an induction combined with \cite[Lemma A-18]{Kumar}.}, Perrin and Smirnov were able to give a general formula for the canonical sheaf of $\overset{\sim}{X_i}$ involving a certain line bundle on $Y$, see \cite[Lemma 4.7]{PerrinSmirnov}.

\begin{prop}\label{prop:canonical sheaf formula}
Let $\mathcal{L}$ be a line bundle on $Y$ such that $q_{i+1}^{*}\L^{-1}$ has degree one on the fibers of $\phi_{i}$. Then \begin{equation}\label{eq:formule generale canonique}
\w_{\overset{\sim}{X}_i}=\O_{\overset{\sim}{X}_i}(-\partial \overset{\sim}{X}_i)\otimes q_i^{*}\L\otimes \Phi_i^{*}(f_0^{*}\L^{-1}\otimes \w_{{\overset{\sim}{X}_0}}).
\end{equation}
\end{prop}

In our specific setting, we have $\w_{\overset{\sim}{X}_0}=\w_{H/B_H}=\L_{H}(-2\rho_H)$. Moreover, we remark that $\Phi_n=\overset{\sim}{p}$ and $f_n=\iota\overset{\sim}{\pi}$ where $\iota$ denotes the immersion $\overline{HB\cdot wB}\hookrightarrow G/B$. Taking $\L=\L_G(-\rho_G)$ in the previous proposition we deduce with $(iv)'$ 
\begin{equation}\label{eq:canonical sheaf expression}
\w_{\overset{\sim}{X}_n}=\O_{\overset{\sim}{X}_n}(-\partial \overset{\sim}{X}_n)\otimes \overset{\sim}{\pi}^{*} \iota^*(\L_{G}(-\rho_G)\otimes \mathcal{M}).
\end{equation}
Thanks to the $B$-canonical splitting of the Bott-Samelson varieties (see \cite[proposition 4.1.17]{BrionKumar}), and by hypotheses $(ii)$ and $(iii)'$, we can still apply \cite[Theorem 20]{HeThomsen} to obtain a splitting of $\overset{\sim}{X_n}$ compatibly splitting $\partial \overset{\sim}{X}_n$. We can then apply the last arguments of Perrin and Smirnov, inspired once again by Brion and Kumar (see \cite[Lemma 5.6]{PerrinSmirnov} with \cite[Theorem 1.2.12]{BrionKumar}), to get the vanishings
$$R^j\overset{\sim}{\pi}_*\O_{\overset{\sim}{X}_n}(-\partial \overset{\sim}{X}_{n})=0,\ j>0,$$ 
and finally the desired (\ref{eq:vanishing canonical sheaf}) by the projection formula on (\ref{eq:canonical sheaf expression}).

\section{Conclusion: Normality, rationality, Cohen-Macaulayness, types A, B, D}\label{sec:conclusion}

In order to obtain Theorem \ref{theo:0 normality, Cohen-Macaulay etc.} and the conclusion in Theorem \ref{theo:1 birational rational} regarding rationality, it is now sufficient to apply Theorem \ref{theo:theorem 2 plus} to the data of $G$, $H$, $B$, $w$ defined by the matrix models and Section \ref{sec:application through matrix models}. Since the birationality assumption $(i)$ is satisfied by the Section \ref{sec:birationality}, it suffices to verify the hypothesis $(ii)$ -or the hypothesis $(iii)'$ in the case of positive characteristic- and the hypothesis $(iv)'$ for the type A. This will be done thanks to the following two propositions. For $i\in \set{1,...,n-1}$ let $\e_i$ denote the morphism 
$$\begin{pmatrix}
t_1 & & \\
& \ddots & \\
& & t_n
\end{pmatrix}\mapsto t_i$$ from $T_n$ to $\mathbb{G}_m$.

\begin{prop}\label{prop:type A caractère dominant}
The character $2\rho_H-\rho_{G\vert T_H}$ is dominant. It is even zero in type A.

\begin{proof}
We begin with the type A case. We have $H\simeq G_r$ and easily 
\begin{align*}
2\rho_H&=\sum_{1\leq i<j\leq r} \e'_i-\e'_j,\\
2\rho_{G\vert T_H}&=\sum_{1\leq i<j\leq n} \e'_i-\e'_j,\\
\e'_i&=\e'_{n-r+i}\ \forall i\in [1,r],\\
\e'_j&=0\ \forall j\in [r+1,n-r].
\end{align*} We deduce
\begin{align*}
2\rho_{G\vert T_H}&=
\sum_{1\leq i<j\leq r} \e'_i-\e'_j\ 
+ \sum_{n-r+1\leq i<j\leq n} \e'_i-\e'_j\\
&+ \sum_{1\leq i\leq r <j\leq n-r} \e'_i-\e'_j\ 
+ \sum_{r+1\leq i\leq n-r<j\leq n} \e'_i-\e'_j \\
&+ \sum_{1\leq i\leq r<n-r+1\leq j\leq n} \e'_i-\e'_j\ +\sum_{r+1\leq i<j\leq n-r} \e'_i-\e'_j \\
&=\sum_{1\leq i<j\leq r} \e'_i-\e'_j\ 
+ \sum_{1\leq l<k\leq r} \e'_{n-r+l}-\e'_{n-r+k}\\
&+ \sum_{1\leq i\leq r <j\leq n-r} \e'_i-0\ 
+ \sum_{1\leq l\leq r<i\leq n-r} 0-\e'_{n-r+l} \\
&+\sum_{1\leq i,l\leq r} \e'_i-\e'_{n-r+l}\ 
+ \sum_{r+1\leq i< j\leq n-r} 0-0\\
&=4\rho_H,
\end{align*}
Hence $2\rho_H-\rho_{G\vert T_H}=0$.
For types B, D, the form of $H$ leads to
\begin{align*}
2\rho_H&=\sum_{1\leq i\leq r/2} (r-2i+2)\e'_i,\\
2\rho_{G\vert T_H}&=\sum_{1\leq i\leq \lfloor n/2\rfloor}(n-2i) \e'_i,\\
\e'_i&=-\e'_{r-i+1}\ \forall i\in [1,r],\\
\e'_j&=0\ \forall j\in [r+1,n-r].
\end{align*} We deduce
\begin{align*}
2\rho_{G\vert T_H}&=\sum_{1\leq i\leq \lfloor n/2\rfloor}(n-2i) \e'_i\\
&=\sum_{1\leq i\leq r}(n-2i) \e'_i\\
&=\sum_{1\leq i\leq r/2}(n-2i)\e'_i+\sum_{r/2+1\leq i\leq r}(n-2i)\e'_i\\
&=\sum_{1\leq i\leq r/2}(n-2i)\e'_i+\sum_{1\leq j\leq r/2}(n-2(r-j+1))\e'_{r-j+1}\\
&=\sum_{1\leq i\leq r/2}\big{(}(n-2i)-(n-2r+2i-2)\big{)}\e'_i\\
&=2\sum_{1\leq i\leq r/2}(r-2i+1)\e'_i\\
&=2\rho_H -\sum_{1\leq i\leq r/2}\e'_i,
\end{align*}
Hence $2\rho_H-\rho_{G\vert T_H}=\sum_{1\leq i\leq r/2}\e'_i$.
\end{proof}
\end{prop}

\begin{rqe}\label{rqe:type C caractère non dominant}
Using similar calculations, we find for type C that $2\rho_H-\rho_{\vert T_H}=-\sum_{1\leq i\leq r/2}\e'_i$, a non-dominant character.
\end{rqe}

\begin{prop}\label{prop:type A surjectivité nabla}
The restriction morphism $$V_G\left((p-1)\rho_{G}\right)\to V_{H}\left((p-1){\rho_{G}}_{\vert T_H}\right)$$ is surjective.
\begin{proof}
We begin again with the type A case. Let us introduce the parabolic subgroup 
$$P'=\Set{ \begin{pmatrix}
A & * & *\\
0 & C  & * \\
0 & 0 & D
\end{pmatrix} \in G| A, D\in G_{r}(R),\ C\in B_{n-2r} }$$ containing $B$ and $$L':= \Set{ \begin{pmatrix}
A & 0 & 0\\
0 & C  & 0 \\
0 & 0 & D
\end{pmatrix} \in G| A,D\in G_{r}(R),\ C\in T_{n-2r} }$$ its Levi subgroup related to $T$. Let us consider the subgroup
$$L'':= \Set{ \begin{pmatrix}
A & 0 & 0\\
0 & I_{n-2r}  & 0 \\
0 & 0 & D
\end{pmatrix} \in G| A,D\in G_{r}(R)},$$
for which $B\cap L''$ is a Borel subgroup. Since $P'/B$ is a Schubert variety, the restriction $$V_G\left((p-1)\rho_{G}\right)\to V_{P'}\left((p-1){\rho_{G}}\right)$$ is surjective (see \cite[Proposition 14.15]{Jantzen1987}). Besides, $(P',L')$ is a Donkin pair in the sense of Donkin (\cite[]{Donkin1985}) because $L'\subset P'$ is a Levi subgroup (\cite[]{Mathieu1990}). Thus, by \cite[Remark 18]{vanderKallen2001}, the restriction $$V_{P'}\left((p-1)\rho_{G}\right)\to V_{L'}\left((p-1){\rho_{G}}\right)$$ is surjective.  But the inclusion $L''\subset L'$ clearly induces an identification $L'/B\cap L'\simeq L''/B\cap L''$ so that $$V_{L'}\left((p-1)\rho_{G}\right)\to V_{L''}\left((p-1){\rho_{G}}_{\vert T\cap L''}\right)$$ is an isomorphism. Since $H\hookrightarrow L''$ is a diagonal embedding, $(L'',H)$ is also a Donkin pair (\cite[]{Mathieu1990} again) and $$V_{L''}\left((p-1){\rho_{G}}_{\vert T\cap L''}\right)\to V_{H}\left((p-1){\rho_{G}}_{\vert T_H}\right)$$ is surjective.  By composition, the desired restriction is therefore surjective.

For types B and D, the argument is exactly the same as for type A except that we replace $P'$, $L'$ and $L''$ by
$$P'=\Set{ \begin{pmatrix}
A & * & *\\
0 & C  & * \\
0 & 0 & {}^{\d}A^{-1}
\end{pmatrix} \in G| A\in G_{r}(R),\ C\in B_{\e,n-2r} },$$
$$L':= \Set{ \begin{pmatrix}
A & 0 & 0\\
0 & C  & 0 \\
0 & 0 & {}^{\d}A^{-1}
\end{pmatrix} \in G| A\in G_{r}(R),\ C\in T_{\e,n-2r} },$$
$$L'':= \Set{ \begin{pmatrix}
A & 0 & 0\\
0 & I_{n-2r}  & 0 \\
0 & 0 & {}^{\d}A^{-1}
\end{pmatrix} \in G| A\in G_{r}(R)},$$
and the pair $(L'',H)$ is Donkin for a different reason. Here, it is no longer because we have a diagonal embedding, but because $H$ is the fixed point set of a diagram automorphism of $L''$ (we get a pair of type $(A_{n-1},C_n)$), see \cite[]{Brundan}.
\end{proof}
\end{prop}

\section{Appendix A: realizations}

\begin{center}\textit{Here $\mathds{k}$ will denote a field that is not necessarily algebraically closed.}\end{center}

\subsection{Definition of the notion} 

For us, realizing a $\k$-variety $X$ or finding a realization of $X$ will simply mean finding finitely presented schemes on relatively large bases whose geometric fibres will provide different incarnations of $X$. The various schemes obtained will be gathered into coherent families. Our point of view is thus a special case of the one adopted in \cite[§8]{EGAIV8a15}, concerning projective families of schemes. As we are dealing with algebras, the term "large" will refer to the order for the inclusion property.

\begin{defi}\label{defi:realization}
Let $X$ be a $\k$-scheme of finite type. A family $(X_A)_{A\in \mathcal{A}}$ is said to be a realization of $X$ if:
\begin{enumerate}
    \item There exists $A_0\in {\mathcal{C}}:=\set{\ZZ\mbox{-subalgebras of }\k \mbox{ of finite type}}$ such that $\mathcal{A}=\set{A\in {\mathcal{C}}|A_0\subset A}$.
    \item For any $A\subset A'$ in $\mathcal{A}$, $X_A$ is a finitely presented $A$-scheme and
    $$X_A\times_A \mathds{k}\simeq X \quad \mathrm{ and }\quad X_A\times_A A' \simeq X_{A'}.$$
\end{enumerate}
\end{defi}

\begin{exemple} Let $X=\mathrm{Spec}\ \mathds{k}[t_1,...,t_n]/(f_1,...,f_r)$. Let $\mathcal{A}$ be the set of all $\ZZ$-subalgebras of $\mathds{k}$ of finite type containing $f_1,...,f_r$. For $A\in \mathcal{A}$, let $X_A:=\mathrm{Spec}\ A[t_1,...,t_n]/(f_1,...,f_r)$.  Then $(X_A)_{A\in \mathcal{A}}$ is a realization of $X$. 
\end{exemple}

The definition \ref{defi:realization} for $\k$-schemes of finite type extends naturally to the notion for morphisms between such schemes, for sheaves and modules, for $\k$-algebraic groups, etc, and for morphisms between these different kinds of objects. In addition, we can require that the objects of the family that \textit{realizes} the data satisfy additional properties, such as being flat on the basis, being a closed immersion, being locally free, etc. Several results ensure the existence of pure realizations for schemes, morphisms, modules, morphisms of modules \cite[Theorem 8.8.2, 8.5.2]{EGAIV8a15} and the fact that we can find several properties for the objects of the family as soon as we consider $A$ sufficiently large \cite[Theorem 8.10.5, Proposition 8.5.5]{EGAIV8a15}. Moreover, enlarging $A$ preserves the commutativity of the diagrams and thus the various algebraic structures involved; in particular, the group structure \cite[Scholia 8.8.3]{EGAIV8a15}. 

For reasons of convenience, we often omit to specify the index set $\mathcal{A}$.

\subsection{Realization of scheme-theoretic images}

We refer to \cite[Section 29.6]{Stacks} for the notion of scheme-theoretic image of a morphism of schemes and its fundamental properties. In particular, we will use the description of such an image when the source of the morphism is reduced \cite[Lemma 29.6.7]{Stacks}. Here we aim to prove Theorem \ref{theo:relèvement image} below and its Corollary \ref{cor:réalisation H1H2...HnZ barre} which provide good realizations for suitable scheme-theoretic images. Our work is motivated by the fact that the construction of such images does not in general commute to (non-flat) base change, so their realization requires special attention. We begin with the following two preliminary propositions.

\begin{prop}\label{prop:dominant image plat intégrité}
Let $X\to S$ be a flat morphism and let $T\to S$ be a schematically dominant quasi-compact morphism.  If $X_T$ is integral, then $X$ is also integral.
\begin{proof}
Consider a cartesian diagram
\[\xymatrix{X_T\ar[r] \ar[d]^{g} & T \ar[d]^f\\
X\ar[r]^u & S, }\]
where $u$ is flat and $f$ is schematically dominant quasi-compact. By flatness, $X\simeq (Im\ f)_X\simeq Im\ g$ and $g$ is also schematically dominant. But, if $X_T$ is integral, then $Im\ g$ is isomorphic to $\overline{g(X_T)}_{red}$ and is also integral. 
\end{proof}
\end{prop}

\begin{prop}\label{prop:formation image ouvert intégrité}
Let $S$ be a scheme and let $f:X\to Y$ be a finitely-presented $S$-morphism with $Im\ f\to S$ open. Then, there exists a nonempty open subset $U\subset S$ such that for any base change $\emptyset\neq T\to U$,  we have
$$ \mathrm{If\ } (Im\ f)_T \mathrm{\ is\ integral\ and\ } X_T \mathrm{\ is\  reduced,\ then\ } Im\ f_T\simeq (Im\ f)_T \mathrm{\ over\ } Y_T. $$
\begin{proof}
Let $\overset{\sim}{f}:X\to Im\ f$ denote the restriction of $f$ onto its scheme-theoretic image. By cancellation, $\overset{\sim}{f}$ is finitely presented just as $f$. Since $f$ is quasi-compact, $Im\ f$ is topologically $\overline{f(X)}=\overline{\overset{\sim}{f}(X)}$ and $\overset{\sim}{f}$ is (topologically) dominant. Since $Im\ f\to S$ is open, we can apply Lemma \ref{lem:realisation image intermediaire} below. It gives a nonempty open subset $U\subset S$ such that for any base change $\emptyset\neq T\to U$, $\overset{\sim}{f}_T(X_T)$ contains a nonempty open subset of $(Im\ f)_T$. If this last scheme is integral then $(\overline{\overset{\sim}{f}_T(X_T)})_{red}\simeq (Im\ f)_T$ and if $X_T$ is reduced then $Im\ \overset{\sim}{f}_T\simeq (\overline{\overset{\sim}{f}_T(X_T)})_{red}$. Under those assumptions $\overset{\sim}{f}_T:X_T\to (Im\ f)_T$ is thus schematically dominant. Since $f_T$ factors as this morphism followed by the closed immersion $(Im\ f)_T\hookrightarrow Y_T$, we deduce that $Im\ f_T\simeq (Im\ f)_T$.
\end{proof}
\end{prop}

\begin{lem}\label{lem:realisation image intermediaire} Let $S$ be a scheme and let $f:X\to Y$ be a finitely presented (topologically) dominant $S$-morphism with $Y\to S$ open. Then, there exists a nonempty open subset $U\subset S$ such that for any base change $\emptyset\neq T\to U$, $f_T(X_T)$ contains a nonempty open subset of $Y_T$. 
\begin{proof}
By Chevalley theorem $f(X)$ is a constructible subset of $Y$ and thus contains a dense open subset $V$ of its closure, which is $Y$ since $f$ is dominant. Since $Y\to S$ is open, $V$ is sent onto a nonempty open subset $U$ of $S$. By restricting over $U$, we get a diagram where the upper horizontal arrows are surjective and the right vertical one is an open immersion:
\[\xymatrix{
f_U^{-1}(V)\ar@{->>}[r] \ar[d] & V \ar@{^{(}->}[d] \ar@{->>}[r] & U\\
X_U\ar[r]^{f_U} & Y_U \ar[ru]}.\]
For any base change $\emptyset\neq T\to U$ the diagram
 \[\xymatrix{
(f_U^{-1}(V))_T\ar@{->>}[r] \ar[d] & V_T \ar@{^{(}->}[d] \ar@{->>}[r] & T\\
X_T\ar[r]^{f_T} & Y_T \ar[ru]}\]
will then preserve those properties. We deduce  $V_T\neq \emptyset$ and that $f_T(X_T)$ contains the open image of $V_T$ in $Y_T$.
\end{proof}

\end{lem}

\begin{theo}\label{theo:relèvement image}
Let $f\colon X\to Y$ be a morphism of $\mathds{k}$-schemes of finite type. Let $(f_A\colon X_A\to Y_A)_A$ be a realization of $f$. Assume that $Y$ is proper over $\mathds{k}$ and that there exists $A$ such that $X_A$ has integral geometric fibers (in particular $X\simeq ((X_A)_{\overline{\eta}})_{\mathds{k}}$ is integral, where $\eta\in \mathrm{Spec}\ A$ denotes the generic point). Then, there exists $A_0$ such that
\begin{enumerate}
\item[(i)]$(Im\ f_A)_{A_0\subset A}$ realizes $Im\ f$.
\item[(ii)]For any $A$ containing $A_0$, we have
\begin{enumerate}
\item $X_A$ is integral and flat over $A$,
\item $Im\ f_A$ is integral with integral geometric fibers and flat, proper, of finite presentation over $A$.
\item The construction of the scheme-theoretic image of $f_A$ commutes with any base change $\emptyset \neq T \to \mathrm{Spec}\ A$ such that $(Im\ f_A)_T$ and $(X_A)_T$ are integral. In particular, for any $s\in \mathrm{Spec} A$, we have an isomorphism over $(Y_A)_{\overline{s}}$:
$$Im\ (f_A)_{\overline{s}}\simeq (Im\ f_A)_{\overline{s}}.$$
\end{enumerate}
\end{enumerate}
\begin{proof}
Let us first show that if $A$ satisfies $(ii)$, then this will also be the case for any $B$ containing $A$ and $(Im\ f_B)_{A\subset B}$ will realize $Im\ f$. By flatness, we have $(Im\ f_A)_\mathds{k}\simeq Im\ (f_A)_\mathds{k}$ and therefore $(Im\ f_A)_\mathds{k}\simeq Im\ f$. Note also that the properties of flatness, properness, being of finite presentation and having integral geometric fibers are preserved for $X_B$ and $(Im\ f_A)_B$ on $B$. These two schemes are also integral by Proposition \ref{prop:dominant image plat intégrité} applied on the two following cartesian diagrams whose right vertical arrows are schematically dominant \[{\xymatrix{X\ar[r] \ar[d] & \mathrm{Spec}\ \mathds{k} \ar[d]\\
X_B\ar[r] & \mathrm{Spec}\ B} 
\quad \xymatrix{Im\ f\ar[r] \ar[d] & \mathrm{Spec}\ \mathds{k} \ar[d]\\
(Im\ f_A)_B \ar[r] & \mathrm{Spec}\ B.}}\] Besides, any base change $\emptyset\neq T\to \mathrm{Spec}\ B$ leads to isomorphisms $$(X_B)_T\simeq ((X_A)_B)_T\simeq (X_A)_T,\ ((Im\ f_A)_B)_T\simeq (Im\ f_A)_T,\ Im\ (f_B)_T\simeq Im\ (f_A)_T.$$ By now applying $(c)$ for $f_A$ and the base change $\emptyset\neq \mathrm{Spec}\ B\to \mathrm{Spec}\ A$ we obtain the desired statements. 

Let us now show the existence of an $A$ satisfying $(ii)$. Since $Y\to \mathrm{Spec}\ \mathds{k}$ is proper, we can assume that $Y_A\to \mathrm{Spec}\ A$ and then $Im\ f_A\to \mathrm{Spec}\ A$ are proper for $A$ sufficiently large (\cite[Theorem 8.10.5]{EGAIV8a15}). We will establish the other assertions by successive localizations. 

\begin{enumerate}
\item[(flatness)] By generic flatness (\cite[Proposition 8.9.4]{EGAIV8a15}) and localization of $A$ at a suitable element,  $Im\ f_A$ and $X_A$ can be supposed to be flat over $A$.
\item[(fiber integrality)] By flatness $(Im\ f_A)_{\overline{\eta}}\simeq Im\ (f_A)_{\overline{\eta}}$. But $Im\ (f_A)_{\overline{\eta}}$ is integral since $(X_A)_{\overline{\eta}}$ is.  Thus the set of $s\in \mathrm{Spec}\ A$ with $(Im\ f_A)_{\overline{s}}$ integral is nonempty. It is an open subset (\cite[Theorem 12.2.4]{EGAIV8a15}) and, localizing again, we can assume that it is the whole of $\mathrm{Spec}\ A$.
\item[(base change)] Since $Im\ f_A\to \mathrm{Spec}\ A$ is open by flatness,  we can apply Proposition \ref{prop:formation image ouvert intégrité} on the finitely presented morphism $f_A:X_A\to Y_A$ and localizing, we can suppose that $(c)$ is satisfied for $f_A$.
\item[(integrality)] We get that $Im\ f_A$ and  $X_A$ are integral by applying Proposition \ref{prop:dominant image plat intégrité} on the following cartesian diagrams whose right vertical arrows are schematically dominant
\begin{equation*}
{\xymatrix{X\ar[r] \ar[d] & \mathrm{Spec}\ \mathds{k} \ar[d]\\
X_A\ar[r] & \mathrm{Spec}\ A} 
\quad \xymatrix{Im\ f\ar[r] \ar[d] & \mathrm{Spec}\ \mathds{k} \ar[d]\\
Im\ f_A \ar[r] & \mathrm{Spec}\ A.}}
\end{equation*}
\end{enumerate}
\end{proof}
\end{theo}

\begin{cor}\label{cor:réalisation H1H2...HnZ barre}
Let $G$ be a $\mathds{k}$-group acting on a $\mathds{k}$-scheme of finite type $X$. Let $Z$ be a closed subscheme of $X$ and $H_1$, $H_2$, ..., $H_n$ be closed subgroups of $G$. Let $(G_A)_A$, $(X_A)_A$, $(G_A\times_A X_A\to X_A)_A$, $(Z_A\hookrightarrow X_A)_A$, $(H_{iA}\hookrightarrow G_A)_A$, $i=1,...n$ be families which realize these data. Assume that for large enough $A$, $Z_A$ and all the $H_{iA}$ have integral geometric fibers and that $X$ is proper over $\mathds{k}$. Then, there exists a family $({Y_A}\hookrightarrow X_A)_A$ which realizes the closed subscheme $(\overline{H_1H_2...H_n\cdot Z})_{red}\hookrightarrow X$ and such that, for any large enough $A$: 
\begin{enumerate}
\item[(i)] $Y_A$ is proper and flat over $A$, with integral geometric fibers.
\item[(ii)] For all $s\in \mathrm{Spec}\ A$, $$({Y_A})_{\overline{s}}\simeq (\overline{ {H}_{1A\overline{s}}{H}_{2A\overline{s}}...{H}_{nA\overline{s}}\cdot Z_{A\overline{s}}})_{red}$$ over ${X}_{A\overline{s}}$.
\end{enumerate}
\begin{proof}
By induction, it suffices to show the result for a single subgroup $H=H_1$ ($n=1$). For large enough $A$ we consider the morphism $f_A\colon H_A\times_A Z_A\to X_A$ which denotes the restriction of the action $G_A\times_A X_A\to X_A$ over $H_A\times_A Z_A\hookrightarrow G_A\times_A X_A$. The family $(f_A)_A$ is clearly a realization of the $\mathds{k}$-morphism of finite type $f:H\times_\mathds{k} Z\to X$ defined in  a similar way by action and restriction. Since the product of two integral schemes over a perfect field is again integral, $H_A\times_A Z_A$ has integral geometric fibers for $A$ sufficiently large. We can therefore apply the previous Theorem \ref{theo:relèvement image} and we get a family of closed immersions $Y_A\hookrightarrow X_A$ realizing $Im\ f\hookrightarrow X$,  with $Y_A$ flat and proper over $A$ and isomorphims $(Y_A)_{\overline{s}}\simeq Im\ (f_A)_{\overline{s}}$ over $(X_A)_{\overline{s}}$ for all $s\in \mathrm{Spec}\ A$. To conclude, all we need to do is to notice that we have $Im\ f\simeq (\overline{H\cdot Z})_{red}$ and $Im\ (f_A)_{\overline{s}}\simeq (\overline{{H_A}_{\overline{s}}\cdot {Z_A}_{\overline{s}}})_{red}$ for all $s\in \mathrm{Spec}\ A$. This comes from the reductiveness of $H_A\times_A Z_A$ and $H\times_\mathds{k} Z$.
\end{proof}
\end{cor}

\subsection{Proof of Lemma \ref{lem:realisation total}}

We can now prove Lemma \ref{lem:realisation total}. Let us assume that $\k$ is algebraically closed and let $G$, $H$, $T$, $B$, $T_H$, $B_H$, $W$, $w_i$ be as in the general setting of Section \ref{sec:context general types free}.

Let first $\underline{G}$, $\underline{H}$, $\underline{T}$, $\underline{B}$, $\underline{T_H}$,  $\underline{B_H}$, $\underline{W}$ be the group schemes over $\mathbb{Z}$ and let $\underline{w}_i\in \underline{W}(\mathbb{Z})$ be the section such that (see \cite[Exposés XIX to XXVI]{SGA3XIXaXXVI}):
\begin{enumerate}
\item[-] $\underline{G}$ is the semi-simple Chevalley group scheme with $\underline{G}_{\mathds{k}}=G$
\item[-] $\underline{H}$ is the reductive Chevalley group scheme with $\underline{H}_{\mathds{k}}=H$
\item[-] $(\underline{T},\underline{B})$ is the Killing pair of $\underline{G}$ with $\underline{T}_{\mathds{k}}=T$, $\underline{B}_{\mathds{k}}=B$
\item[-] $(\underline{T_H},\underline{B_H})$ is the Killing pair of $\underline{H}$ with $\underline{T_H}_{\mathds{k}}=T_H$, $\underline{B_H}_{\mathds{k}}=B_H$
\item[-] $\underline{W}$ is the Weyl group of $\underline{G}$ related to $\underline{T}$
\item[-] The base change ${\underline{w}_i}_{\mathds{k}}$ is the element $w_i\in W$
\end{enumerate} 
We know that the fppf quotient sheaf $\underline{G}/ \underline{B}$ is representable by a $\mathbb{Z}$-scheme of finite presentation, smooth and projective and such that its base change $(\underline{G}/ \underline{B})_{K}$ over any algebraically closed field $K$ is the flag variety $ \underline{G}_{K}/ \underline{B}_{K}$ (see \cite[Exposé XXIV, Théorème 1.3]{SGA3XIXaXXVI}). 

We now make use of \cite[§8]{EGAIV8a15}. By enlarging $A$ at each step, we successively prove the following assertions. Since there is a closed immersion of groups $H\hookrightarrow G$, we can suppose that $\underline{H}_A$ is a closed subgroup of $\underline{G}_A$. Besides, since the base changes over $\k$ of $\underline{T_H}_A$ and $\underline{T}_A\times_{\underline{G}_A}\underline{H}_A$ give $T_H$, we can assume that these two groups identify over $\underline{H}_A$. Similarly $\underline{B_H}_A$ and $\underline{B}_A\times_{\underline{G}_A}\underline{H}_A$ can be identified over $\underline{H}_A$. Finally, the root data are preserved and, for any $s\in \mathrm{Spec}\ A$, $$2\rho_{\underline{H}_{A\overline{s}}}-{\rho_{\underline{G}_{A\overline{s}}}}_{\vert \underline{T_H}_{A\overline{s}}} \mbox{ is dominant.}$$
On the other hand, there exist semi-ample line bundles $\mathcal{L}_A$ on $(\underline{G}/\underline{B})_A$, for $A$ sufficiently large so that
$$ (\mathcal{L}_A)_A \mbox{ realizes } \mathcal{L}.$$ 

To finish the proof of the lemma, we now just need to find a good realization of our variety $X_i=\overline{HB\cdot w_iB}$. This will be done thanks to the following consequence of Corollary \ref{cor:réalisation H1H2...HnZ barre}. 

\begin{cor}\label{cor:réalisation HBwB barre}
There exists $\left({Y_{iA}}\hookrightarrow (\underline{G}/\underline{B})_A\right)_A$ which realizes the closed subvariety $X_i\hookrightarrow G/B$ and such that, for large enough $A$
\begin{enumerate}
\item[(i)] $Y_{iA}$ is projective and flat over $A$,

\item[(ii)] For all $s\in \mathrm{Spec}\ A$, $$({Y_{iA}})_{\overline{s}}\simeq \overline{ {\underline{H}}_{A\overline{s}}{\underline{B}}_{A\overline{s}}\cdot \underline{w}_{A\overline{s}}{\underline{B}}_{A\overline{s}}}$$  over $\underline{G}_{A\overline{s}}/{\underline{B}}_{A\overline{s}}$.
\end{enumerate}
\begin{proof}
Suppose $A$ is large enough to guarantee all the above statements and existences. Let $Z:=\mathrm{Spec}\ \k$ and $Z\hookrightarrow G/B$ be the section which corresponds to the closed point $w_i B\in G/B$. The choice of $\underline{B}_A$ gives a morphism $\underline{W}_A\to (\underline{G}/\underline{B})_A$ defined by the natural transformation $(n\underline{T}_A(S)\mapsto n\underline{B}_A(S))_S$ on the corresponding sheaves. Composing with  $\underline{w}_{iA}\colon \mathrm{Spec}\ A\to \underline{W}_A$ produces a section $\mathrm{Spec}\ A\to (\underline{G}/\underline{B})_A$.  Its geometric fiber over $s\in \mathrm{Spec}\ A$ is integral and corresponds to the closed point ${\underline{w}_{iA}}_{\overline{s}}\underline{B}_{A{\overline{s}}}\in (\underline{G}/\underline{B})_{A{\overline{s}}}$. Besides, the family of these sections realizes $Z\hookrightarrow G/B$. Moreover $(\underline{G}/\underline{B})_A$ is proper over $A$ and $\underline{B}_A$ and $\underline{H}_A$ have integral geometric fibers as Borel and reductive group schemes. Finally, the natural action of $G$ on $G/B$ is realized by the family of the natural actions of $\underline{G}_A$ on $(\underline{G}/\underline{B})_A$. By applying Corollary \ref{cor:réalisation H1H2...HnZ barre} with $H_1:=H$, $H_2:=B$ we obtain all the desired assertions, except for the projectivity of the $Y_{iA}$ which is automatically satisfied as they are closed subschemes of $(\underline{G}/\underline{B})_A$. 
\end{proof}
\end{cor}

\section{Appendix B: on Perrin and Smirnov's arguments}

\begin{center}\textit{Here $\k$ denotes an algebraically closed field with }$char(\k)\neq 2$.\end{center}

\subsection{An identical birational morphism}

We present Perrin and Smirnov's construction of a birational morphism onto an irreducible component of a Springer fiber, in type A. We also show that these morphisms are special cases of those we have constructed.

\paragraph{Perrin and Smirnov's version}

Let $V$ be a $n$-dimensional $\k$-vector space. Let $N$ be a nilpotent endomorphism of $V$ of nilpotency order two, with rank $r$ and let $Z_N$ be its centralizer in the general linear group $Gl(V)$. In this setting, let 
$$\mathcal{F}:=\mathcal{F}(V)=\Set{V_1\subset \dots \subset V_{n-1}| V_i\subset V,\ \dim V_i=i,\ \forall
 i}$$ be the flag variety and let $$\mathcal{F}_N:=\Set{V_{\bullet}\in \mathcal{F}| N(V_i)\subset V_i,\ \forall i}$$ be the Springer fiber over $N$.  Let $$\tau=\begin{ytableau} 
n & p_r \\ 
\none[\vdots] & \none[\vdots] \\ 
\cdot & p_1 \\ 
\none[\vdots] & \none\\
1
\end{ytableau}$$
be a standard tableau (with decreasing numbers from left to right and from top to bottom) and let
$$X:=X_{\tau}=\Set{ F_{\bullet}\in \F (V)| \dim F_{p_i}\cap Im\ N\geq i,\ \dim F_{p_i}\cap Ker\ N\geq p_i-i+1,\ F_{p_i}\subset N^{-1}(F_{p_{i-1}}), \forall i\in [1,r]}$$ be the related irreducible component of the Springer fiber over $N$. We define then the variety
$$\hat{X}:=\hat{X}_{{\tau}}=\big{\lbrace}(F'_{\bullet},F_{\bullet})\in \F (Im\ N)\times \F (V)\vert \ F'_i\subset F_{p_i}\subset N^{-1}(F'_{i-1}) \ \forall i\in [1,r] \big{\rbrace}.$$

In \cite{PerrinSmirnov}, Perrin and Smirnov show that $\hat{X}$ is smooth, irreducible and that the projection to $\F$ induces a proper birational $Z_N$-equivariant morphism 
\begin{equation}\label{eq:PerrinSmirnov morphism}
\hat{X}\to X
\end{equation}
as soon as $p_{i+1}>p_i+1$ for all $i$ (if this is not the case, the result is valid for an irreducible component of $\hat{X}$). 

\paragraph{Our version} Let $G$, $H$, $B_H$, $Z$, $e$ be as in the matrix setting of Section \ref{sec:birationality} (type A). Let be integers $q_1,...,q_r$ and $s_1,...,s_{n-2r}$ such that
\begin{enumerate}
\item[-] $q_r:=p_r+1=min \big{(} p_r+1,...,n \big{)} $ and by decreasing induction on $i=r-1,...,1$, $q_i:=min\big{(} \lbrace p_i+1,...,n\rbrace \setminus \lbrace p_{i+1},q_{i+1},...,p_r,q_r\rbrace \big{)} $.
\item[-] $\lbrace s_1,...,s_{n-2r}\rbrace:=\lbrace 1,...,n\rbrace \setminus \lbrace p_{1},q_{1},...,p_r,q_r \rbrace$ with $s_i<s_{i+1}$.
\end{enumerate}
We then define $w:=w_{\tau}$ the unique element of $\mathfrak{S}_n$ such that
\begin{enumerate}
\item[-] $w(p_i):=i$, $\forall i\in \lbrace 1,...,r\rbrace $.
\item[-] $w(q_i):=n-r+i$, $\forall i\in \lbrace 1,...,r\rbrace$.
\item[-] $w(s_j):=n-r+1-j$, $\forall j\in \lbrace 1,...,n-2r\rbrace$.
\end{enumerate}
Let us remark that $w$ induces the increasing bijection from $\lbrace p_1,...,p_r\rbrace$ to $\lbrace 1,...,r\rbrace $ and the decreasing bijection from $\lbrace s_1,...,s_{n-2r}\rbrace $ to $\lbrace r+1,...,n-r\rbrace$. Here are some examples:
\begin{align*}
a)\ &\tau=\begin{ytableau} 
5 & 4 \\ 
 3 & 2 \\ 
1 \\ 
\end{ytableau} &w_{\tau}=
\begin{pmatrix}
1 & 2 & 3 & 4 & 5 \\
3 & 1 & 4 & 2 & 5 
\end{pmatrix}.\\
b)\  &\tau=\begin{ytableau} 
6 & 4 \\ 
 5 & 2 \\ 
3 \\
1 
\end{ytableau} &w_{\tau}=
\begin{pmatrix}
1 & 2 & 3 & 4 & 5 & 6 \\
4 & 1 & 5 & 2 & 6 & 3 
\end{pmatrix}.\\
c)\ &\tau=\begin{ytableau} 
7 & 5 \\ 
 6 & 4 \\ 
3 & 2 \\
1 
\end{ytableau} &w_{\tau}=
\begin{pmatrix}
1 & 2 & 3 & 4 & 5 & 6 & 7\\
4 & 1 & 5 & 2 & 3 & 7 & 6
\end{pmatrix}.
\end{align*}

We can check that $w$ satisfies the conditions of Proposition \ref{prop:type A a b quasi imply hypo 1} so that we can apply the proof of Section \ref{sec:birationality} and get a birational $Z$-equivariant morphism as (\ref{eq:morphisme birationnel}): 
\begin{equation}\label{eq:morphism birationnel type A special w}
H\times ^{B_H}\overline{B\cdot wB}\to \overline{Z\cdot wB}.
\end{equation}

Let us show that the morphisms (\ref{eq:PerrinSmirnov morphism}) and (\ref{eq:morphism birationnel type A special w}) identify. 

\paragraph{Identification} 
There is a basis $(f_i)$ of $V$ such that:
\begin{align*}
ImN&=\langle f_1,...,f_{r}\rangle\\
KerN&=\langle f_1,...,f_{n-r}\rangle\\
N(f_{n-r+i})&=f_i, \ \forall i\in \set{1,...r}.
\end{align*}
By choosing it so that $G$ acts on $V$ and then on $\mathcal{F}$, we identify $G$ with $Gl(V)$, $e$ with $N$, $Z$ with $Z_N$, $H$ with a diagonal embedding of $Gl(Im\ N)$ into $Gl(V)$ and $B$ with the stabilizer of $F:= \langle f_1 \rangle\subset \dots \subset \langle f_1,...,f_{n-1}\rangle \in \mathcal{F}$. Hence, in particular, $$G/B\simeq \mathcal{F}.$$ Besides, the element $w$ act on $\F$. Note that
$$wF\in X$$ so that $Z\cdot wF\subset X$. Since $$\dim Z\cdot wB = \dim H/B_H+\ell (w)$$ thanks to (\ref{eq:morphism birationnel type A special w}), the computations of $\dim H/B_H=\binom{r}{2}$ and $\ell (w)=\binom{n-r}{2}$ lead to the equality $$\overline{Z\cdot wF}= X.$$
On the other hand, we can check (when $p_{i+1}>p_i+1$ for all $i$) that the fiber of (\ref{eq:PerrinSmirnov morphism}) over $wF$ is isomorphic to the Schubert variety defined by $w$, i.e. $$\overline{B\cdot wF}.$$ This gives the following diagram \begin{equation}\label{eq:identification birational Perrin Smirnov}
\xymatrix{H\times^{B_H} \overline{B\cdot wB}\ar[d]_{\wr} \ar[r] & \overline{Z\cdot wB}\ar[d]_{\wr}  \ar@{^{(}->}[r] &  G/B\ar[d]_{\wr} \\
\hat{X}\ar[r] & X\ar@{^{(}->}[r] & \mathcal{F},}\end{equation}
identifying (\ref{eq:morphism birationnel type A special w}) and (\ref{eq:PerrinSmirnov morphism}) through the base $(f_i)$ and through a natural isomorphism between $H$-equivariant bundles with isomorphic fibers. 

\subsection{The problem of normality, a counter-example}

While the general arguments presented in Proposition \ref{prop:PerrinSmirnov induction} and \ref{prop:canonical sheaf formula} are taken from the article of Perrin and Smirnov, the authors actually apply them in a different way for the type D. The reason is that they embed the irreducible component of the Springer fiber into a larger variety, which does not live in the flag variety of the ambient group but in its product with a Lagrangian space (see their Proposition 3.18). They thus deal with a $X_n$ different from ours, which allows them to obtain a formula analogous to (\ref{eq:canonical sheaf expression}) without any additional assumption on the pullback sheaf like hypothesis $(iv)'$ of Theorem \ref{theo:theorem 2 plus} (see their Lemma 4.7). The problem is that the embedding they present does not exist in general as an algebraic morphism, so they cannot carry out their argument and even ensure normality. Let us present a counter-example in the context of their article.

\paragraph{Context}

Let $V$ be a $2n$-dimensional $\k$-vector space, and $SO(\w)$ be the group of unimodular linear operators preserving a symmetric nondegenerate bilinear form $\w$. Let $N\in \mathfrak{g}$ be a nilpotent antiadjoint endomorphism.  Let $Z_N$ be the stabilizer of $N$ in $SO(\w)$: $$Z_N:=\Set{g\in SO(\w)| gNg^{-1}=N}.$$ If $2r=\dim ImN$, there exists a basis $(f_i)$ of $V$ such that:
\begin{align*}
ImN&=\langle f_1,...,f_{2r}\rangle\\
KerN&=\langle f_1,...,f_{2n-2r}\rangle\\
N(f_{2n-2r+i})&=f_i, \ \forall i\in \set{1,...r}\\
N(f_{2n-2r+i})&=-f_i, \ \forall i\in \set{r+1,...,2r}\\
\w(f_i,f_j)&=\delta_{i,2n-j+1},\ \forall i,j\in \set{1,...,2n}
\end{align*}

Let now $\mathcal{OF}$ be the variety of orthogonal flags defined by:

$$\mathcal{OF}:=\Set{V_1\subset ...\subset V_{n-1}\subset V_{n+1}\subset...\subset  V_{2n-1}| V_{2n-i}=V_i^{\perp},\ \dim V_i=i\ \forall i },$$
and $\mathcal{OF}_N$ be the closed subvariety of $N$-stable flags, i.e. the Springer fiber over $N$:

$$\mathcal{OF}_N:=\Set{V_{\bullet}\in \mathcal{OF}| N(V_i)\subset V_i\ \forall i}.$$

We can endow the vector space $ImN$ with a skew-symmetric nondegenerate bilinear form $\a$ satisfying:
$$\alpha(u,N(v))=\w(u,N(v))\ \forall u\in ImN ,\ v\in V.$$

Let then $\mathcal{L}$ be the variety of Lagrangian subspaces of $ImN$ for $\a$:

$$\mathcal{L}:=\set{W\subset ImN| W \mbox{ is totally } \a \mbox{-isotropic}}$$

If $(f'_i)$ is any basis of $V$ such that $\w(f'_i,f'_j)=0$ for $i+j\neq 2n+1$, we will denote by $F(f'_1,...,f'_{2n})$ the flag in $\mathcal{OF}$ such that, for all $i$, $$F(f'_1,...,f'_{2n})_i:=\langle f'_1,...,f'_{i}\rangle,$$ and merely by $F_{\bullet}$ the flag $F(f_1,...,f_{2n})$. We consider the application $\phi$ defined as follows: $$\phi:\ \mathcal{OF}_N\to \mathcal{L},\ V_{\bullet}\mapsto \sum_{i=1}^{n-1}V_i\cap N(V_i^{\perp})$$ According to \cite[Remark 3.13]{PerrinSmirnov}, $\phi$ is well-defined. It is clear that $Z_N$ acts on $\mathcal{OF}$ and $\mathcal{L}$ and that $\phi$ is $Z_N$-equivariant.

\paragraph{Non-continuity}

Now, if the embedding of the authors exists as an algebraic morphism, then $\phi$ must exist as well. But we will show that $\phi$ is not continuous in the case $n=2r=4$.

Let $w$ and $s$ be the operators in $SO(\w)$ which correspond respectively to the permutation $15263748$ and $13245768$\footnote{Let $T$ be the maximal torus of $SO(\w)$ related to $(f_i)$,  $\e_i$ be the characters on $T$ defined by $t\mapsto f_i^*(t(f_i))$ and $W$ be the Weyl group attached to $T$. Let $s_2$ and $s_4$ in $W$ be respectively the reflections associated to the roots $\epsilon_2-\epsilon_{3}$ and $\epsilon_3+\epsilon_{4}$. Then, $w$ and $s$ can be respectively seen as representatives of the Weyl group elements $s_4s_2$ and $s_2$.} of $(f_i)$, and let $\set{U(t)}$ be the one-parameter subgroup of $Z_N$ acting on $(f_i)$ with the matrix

$$U(t):=\begin{pmatrix}
1 & & & & & & &\\
& 1 & & & & & &\\
 & & 1 & 0  & t &  0& & \\
 & &  & 1 & 0 & -t & & \\
 & &  &  & 1 & 0 & &\\
  & &  &  &  & 1 &  & \\
  & & & & & & 1 &\\
  & & & & & & & 1
  \end{pmatrix}.$$

We can check that for all $t\neq 0$, 

$$U(t)wF_{\bullet}=F(f_1,tf_3+f_5,f_2,-tf_4+f_6,f_3,f_7,f_4,f_8)=F(f_1,f_3+1/tf_5,f_2,f_4 -1/tf_6,f_3,f_7,f_4,f_8)$$

and  
$$sF_{\bullet}=F(f_1,f_3,f_2,f_4,f_5,f_7,f_6,f_8),$$

and we see that 
$$lim_{t\to \infty}\ U(t)wF_{\bullet}=sF_{\bullet}.$$

But $wF_{\bullet}$ is $N$-stable and
$$\phi(wF_{\bullet})=\langle f_1,f_2\rangle.$$

We have also
$$\phi(U(t)wF_{\bullet})=U(t)\langle f_1,f_2\rangle=\langle f_1,f_2\rangle$$

and 
$$\phi(sF_{\bullet})=s\langle f_1,f_2\rangle=\langle f_1,f_3\rangle.$$

Hence, 
$$lim_{t\to \infty}\ \phi(U(t)wF_{\bullet})= \langle f_1,f_2\rangle \neq \langle f_1,f_3\rangle=\phi(lim_{t\to \infty}\ U(t)wF_{\bullet}).$$

\subsection{Some additional comments}

\paragraph{} In our context, the Lagrangian space of Perrin and Smirnov corresponds exactly to $H/P_H$, where $P_H$ is the parabolic subgroup of $H$ containing $B_H$, and characterized by the $(r/2)^{th}$ simple root ($\e'_{r/2}-\e'_{r/2+1})$ according to Section \ref{sec:conclusion}). The problem of embedding amounts then to the existence of an algebraic morphism $\psi$ from $G/B$ to $H/P_H$ making the following diagram commutative \[\xymatrix{
H\times^{B_H}\overline{B\cdot wB} \ar[d]_-{q} \ar[r]^-{k}& H/B_H\ar@{^{(}->}[d]\\
G/B\ar[r]^-{\psi} & H/P_H.}\] Let us remark that the line bundle $\L_H(\rho_{G\vert T_H}-2\rho_H)$ on $H/B_H$ is the pullback of an equivariant line bundle on $H/P_H$ so that such a diagram would guarantee that we satisfy the hypothesis $(iv)'$ of Theorem \ref{theo:theorem 2 plus}. Therefore the embedding problem seems to be crucial for obtaining rational resolutions and the Cohen-Macaulay property, in Perrin and Smirnov's argument as well as in ours. While the previous counter-example does not rigorously forbid the existence of $\psi$, it does indicate that another approach might be necessary if we want to prove these two additional results.

\printbibliography

\end{document}